\documentclass[12pt, reqno]{amsart}

\usepackage{latexsym, amsfonts, amsmath, amsthm, amscd} 
\newcommand{\tsig}{\widetilde{\sigma}}
\newcommand{\ttau}{\widetilde{\tau}} 

\newcommand{\pair}[1]{\langle#1\rangle}
 
\newcommand{\FD}[2][]{\mathfrak{frD}_{#1}(#2)}
\newcommand{\FDpp}[1][]{\mathfrak{frD}_{#1}(\pp)}
\newcommand{\FDpn}[1][]{\mathfrak{frD}_{#1}(\pn)}

\newcommand{\Atop}[2]{\genfrac{}{}{0pt}{}{#1}{#2}}

\newcommand{\jtp}[1]{\{\mkern-8.5mu\{#1\}\mkern-8.5mu\}}
\newcommand{\FT}{{\mathcal F}}

\newcommand{\Elaood}{(\Elao)\od}\newcommand{\Elaoev}{(\Elao)\ev}

\newcommand{\RModd}{\RM\odd}
       \newcommand{\RMeven}{\RM\even}

\newcommand{\Eodd}{\E\odd}  \newcommand{\Eeven}{\E\even}
\newcommand{\Dodd}{\DD\odd} 
\newcommand{\Deven}{\DD\even}

\newcommand{\Elaod}{(\Ela)\od}
\newcommand{\Elaev}{(\Ela)\ev}

\newcommand{\cHod}{\cH\od} 
\newcommand{\cHev}{\cH\ev}  
\newcommand{\cHp}{\cH'}
\newcommand{\cHpod}{(\cH')\od} 
\newcommand{\cHpev}{(\cH')\ev}   
\newcommand{\ev}{^{\uparrow}} 
\newcommand{\od}{^{\downarrow}} 
\newcommand{\evod}{^{\updownarrow}} 
\newcommand{\even}{^{even}} 
\newcommand{\odd}{^{odd}}

\newcommand{\lao}{\la_0} \newcommand{\laop}{\lao'}
 
\newcommand{\Elao}{\E^{\lao}} \newcommand{\Elaop}{\E^{\laop}}
 
\newcommand{\Blao}{\B^{\lao}} \newcommand{\Blaop}{\B^{\laop}}
\newcommand{\galao}{\ga_{\lao}}
\newcommand{\tgala}{\tga_{\la}}
\newcommand{\tgalao}{\tga_{\lao}}
\newcommand{\gala}{\ga_{\la}}
\newcommand{\Wla}[1][^{\la}]{W^{#1}}

\newcommand{\pilao}{\pi_{\lao}}
\newcommand{\Jlao}{J^{\lao}}\newcommand{\Jlaop}{J^{\laop}}
\newcommand{\bhla}{{\mathbf h}_{\la}}

\newcommand{\gfin}{^{\g-fin}}
\newcommand{\sgfin}{_{\g-fin}}
 
\newcommand{\tsum}{\textstyle\sum}
\newcommand{\tK}{\widetilde{K}}

\newcommand{\boxt}{\Tox\kern -6.3pt\raise .55pt
    \hbox{$\scriptstyle{\times}$}}%boxtimes

\newcommand{\bfb}{{\mathbf b}} 
                 \newcommand{\tbfb}{\widetilde{\bfb}}

\newcommand{\bm}{{\mathbf m}} 
\newcommand{\bp}{{\mathbf p}} 
 \newcommand{\bq}{{\mathbf q}} 
\newcommand{\bj}{{\mathbf j}}

\newcommand{\bsla}{{\mathbf s}^{\la}}
\newcommand{\bvla}{\bv^{\la}} 
\newcommand{\btt}{{\mathbf t}} 

\newcommand{\bv}{{\mathbf v}}
\newcommand{\bx}{{\mathbf x}}

\newcommand{\Innw}{\mathop{\mathrm{Inn}}\nolimits_{w}}

\newcommand{\cdott}{\cdot,\cdot}

\newcommand{\D}{{\mathfrak D}}
\newcommand{\Dpp}[1][]{\D_{#1}(\pp)} 
\newcommand{\Dpn}[1][]{\D_{#1}(\pn)}
\newcommand{\Tpp}{T^*\,{\pp}}  
\newcommand{\Tpn}{T^*\,{\pn}}

\newcommand{\RTpp}{R(T^*\,{\pp})}  
\newcommand{\RTpn}{R(T^*\,{\pn})}
\newcommand{\DtZ}[1][]{\D_{#1}(\tZ)}  
\newcommand{\DZ}[1][]{\D_{#1}(Z)}

\newcommand{\grs}{\gr^{\sharp}}
\newcommand{\bps}{\bp^{\sharp}}
 
\newcommand{\DtZs}[1][]{\D^{\sharp}_{#1}(\tZ)}  
\newcommand{\DZs}[1][]{\D^{\sharp}_{#1}(Z)}  

\newcommand{\Bla}{\B^{\la}} 
\newcommand{\Ela}{\E^{\la}} 
\newcommand{\Cla}{\mathcal C^{\la}}
\newcommand{\RO}[1][]{R^{#1}(\OO)} 
\newcommand{\RM}[1][]{R^{#1}(\M)} 
 
\newcommand{\RTX}[1][]{R^{#1}(T^*X)}
\newcommand{\RZ}[1][]{R^{#1}(Z)}
\newcommand{\RtZ}[1][]{R^{#1}(\tZ)}
\newcommand{\RTZ}[1][]{R^{#1}(T^*Z)}
\newcommand{\RTtZ}[1][]{R^{#1}(T^*\tZ)}
\newcommand{\Spn}[1][]{S^{#1}(\pn)}
\newcommand{\Spp}[1][]{S^{#1}(\pp)} 
\newcommand{\Spm}[1][]{S^{#1}(\p^\pm)}
\newcommand{\Smp}[1][]{S^{#1}(\p^\mp)}

\newcommand{\TZ}{T^*Z}
\newcommand{\TtZ}{T^*\tZ}
 
\newcommand{\Dla}{\D^{\la}}  
\newcommand{\Dlao}{\D^{\lao}}  \newcommand{\Dlaop}{\D^{\laop}}
\newcommand{\phiF}{\phi^{F}} 
\newcommand{\tB}[1][]{\widetilde{\B}_{#1}}

\newcommand{\tom}{\widetilde{\om}} 
\newcommand{\tL}{\widetilde{L}} 
\newcommand{\tZ}{\widetilde{Z}}  
 
\newcommand{\tga}{\widetilde{\ga}}

\newcommand{\Pila}{\Pi_{\la}} \newcommand{\Pilap}{\Pi_{\la'}} 
\newcommand{\Pilao}{\Pi_{\lao}}

\newcommand{\Jla}{J^{\la}} \newcommand{\Jlat}{J^{1-\la}}  
\newcommand{\etala}{\eta_{\la}} 

\newcommand{\etalas}{\etala^{\bullet}} 
\newcommand{\etalats}{\eta_{1-\la}^{\bullet}}
\newcommand{\pila}{\pi_{\la}}
\newcommand{\pilas}{\pila^{\bullet}} 
\newcommand{\mus}{\mu^{\bullet}}

\newcommand{\nr}{m} 
\newcommand{\cla}{c_{\la}}  
\newcommand{\Ila}{\cJ_{\la}}  
\newcommand{\pn}{\p^{-}} 
\newcommand{\pp}{\p^{+}}
\newcommand{\ppn}{\p^{\pm}}
\newcommand{\pnp}{\p^{\mp}}
 
\newcommand{\halfN}{\half{\mathbb N}}
\newcommand{\halfZ}{\half{\mathbb Z}}
\newcommand{\delF}{\del_{F}} 

\newcommand{\ks}{\kk^{\s}}

\newcommand{\TT}{{\mathfrak T}} 
 
\newcommand{\btheta}{\underline{\theta}} 
 
\newcommand{\ltimes}{\vbox to 5.4pt{\leaders\vrule\vfil}\kern
  -5pt\times}

\newcommand{\del}{\partial} 
\newcommand{\pd}[2]{\frac{\partial#1}{\partial#2}}  

\newcommand{\pdb}[2][{}]
          {\frac{\partial^{#1}\phantom{x}}{\partial#2}}

\newcommand{\llongrightarrow}{\relbar\joinrel\longrightarrow}
\newcommand{\lto}{\longrightarrow}  
\newcommand{\hhookrightarrow}{\lhook\joinrel
    \llongrightarrow}
\newcommand{\mapright}[1]{\smash{\mathop
   {\longrightarrow}\limits^{#1}}}
\newcommand{\mmaprightd}[1]{\smash{\mathop
   {\llongrightarrow}\limits^{#1}}}   
\newcommand{\mmapright}[1]{\;\smash{\mathop
   {\llongrightarrow}\limits^{#1}}\;}

\newcommand{\llongleftarrow}{\longleftarrow\joinrel\relbar}
\newcommand{\mmapleft}[1]{\;\smash{\mathop
   {\llongleftarrow}\limits^{#1}}\;}

\newtheorem{thm}{Theorem}[subsection]
\newtheorem{lem}[thm]{Lemma}

\newtheorem{prop}[thm]{Proposition}
\newtheorem{cor}[thm]{Corollary}

\theoremstyle{definition}
\newtheorem{definition}[thm]{Definition}
\newtheorem{example}[thm]{Example}

\theoremstyle{remark}
\newtheorem{remark}[thm]{Remark}

\numberwithin{equation}{section}

\newcommand{\mapdown}[1]{\Big\downarrow\rlap{$\vcenter
{\hbox{$\scriptstyle#1$}}$}}
\newcommand{\mapup}[1]{\Big\uparrow\rlap{$\vcenter
{\hbox{$\scriptstyle#1$}}$}}

\newcommand{\inv}{^{-1}} 
\newcommand{\ovl}{\overline}

\newcommand{\half}[1][1]{\frac{#1}{2}}
\newcommand{\thalf}{\textstyle\frac{1}{2}} 
\newcommand{\four}{\frac{1}{4}} 
\newcommand{\tfour}{\textstyle\frac{1}{4}}

\newcommand{\fa}{{\mathfrak a}}
\newcommand{\fb}{{\mathfrak b}} 
\newcommand{\fe}{{\mathfrak e}}

          \newcommand{\g}{{\mathfrak g}}
\newcommand{\fh}{{\mathfrak h}}
 
\newcommand{\fk}{{\mathfrak k}} 
    \newcommand{\kk}{{\mathfrak k}}

          \newcommand{\p}{{\mathfrak p}}
          \newcommand{\q}{{\mathfrak q}}
\newcommand{\fr}{{\mathfrak r}}
              \newcommand{\s}{{\mathfrak s}}

\newcommand{\fso}{\mathfrak{so}} 
\newcommand{\fsu}{\mathfrak{su}} 
\newcommand{\fsp}{\mathfrak{sp}}
\newcommand{\fsl}{\mathfrak{sl}}
\newcommand{\fgl}{\mathfrak{gl}}

\newcommand{\ove}{{\ovl{e}}}

\newcommand{\C}{{\mathbb C}}
\newcommand{\Z}{{\mathbb Z}}
\newcommand{\R}{{\mathbb R}}
  
\newcommand{\bbN}{{\mathbb N}} \newcommand{\bbT}{{\mathcal T}} 

\newcommand{\gR}{\g_{\R}}
\newcommand{\GR}{G_{\R}}
\newcommand{\KR}{K_{\R}}
\newcommand{\pR}{\p_{\R}}
\newcommand{\kR}{\fk_{\R}} 
\newcommand{\OR}{\OO_{\R}}
 
\newcommand{\om}{\omega}  \newcommand{\Om}{\Omega}
\newcommand{\al}{\alpha}
\newcommand{\be}{\beta}
\newcommand{\ga}{\gamma}

\newcommand{\kap}{\kappa}
\newcommand{\la}{\lambda} \newcommand{\La}{\Lambda}
\newcommand{\vrho}{\varrho}
\newcommand{\sig}{\sigma}  
\newcommand{\vsig}{\varsigma}

\newcommand{\ze}{\zeta}

                                     \newcommand{\A}{{\mathcal A}}
                                     \newcommand{\B}{{\mathcal B}} 
\newcommand{\cC}{{\mathcal C}}

\newcommand{\cG}{{\mathcal G}}
\newcommand{\cH}{{\mathcal H}}
                                             \newcommand{\I}{{\mathcal I}}
\newcommand{\cJ}{{\mathcal J}}
                                         
\newcommand{\cL}{{\mathcal L}}
                                                \newcommand{\M}{{\mathcal M}}
\newcommand{\cN}{{\mathcal N}} 
                                          \newcommand{\OO}{{\mathcal O}} 
\newcommand{\cP}{{\mathcal P}} 
                                            \newcommand{\Q}{{\mathcal Q}} 
\newcommand{\cS}{{\mathcal S}} 
                                         \newcommand{\T}{{\mathcal T}}  
                                           \newcommand{\U}{{\mathcal U}} 
\newcommand{\cV}{{\mathcal V}}

\newcommand{\aND}{\quad\mbox{and}\quad} 
\newcommand{\AND}{\qquad\mbox{and}\qquad} 
\newcommand{\Cl}{\mathrm{Cl}}

\newcommand{\adj}{\mathop{\mathrm{ad}}\nolimits}

\newcommand{\gr}{\mathop{\mathrm{gr}}\nolimits}
\newcommand{\ugr}{\mathop{\underline{\mathrm{gr}}}\nolimits}

\newcommand{\Cent}{\mathop{\mathrm{Cent}}\nolimits}
\newcommand{\End}{\mathop{\mathrm{End}}\nolimits}
\newcommand{\Hom}{\mathop{\mathrm{Hom}}\nolimits}

\newcommand{\rank}{\mathop{\mathrm{rank}}\nolimits}
\newcommand{\Tr}{\mathop{\mathrm{Tr}}\nolimits}
\newcommand{\Det}{\mathop{\mathrm{Det}}\nolimits}
\newcommand{\tr}{\mathop{\mathrm{tr}}\nolimits}

\newcommand{\Vect}{\mathop{\mathfrak{Vect}}\nolimits}

\newcommand{\dd}{\mathord{\mathrm{d}}}   
\newcommand{\gog}{\g\oplus\g}

\newcommand{\Ug}[1][]{{\mathcal U}_{#1}({\mathfrak g})} 
 
\newcommand{\Sg}{S({\mathfrak g})}  
\newcommand{\RR}[1][]{{\mathcal R^{#1}}}

\newcommand{\Killg}[1]{({#1})_{\g}}

\newcommand{\NN}{{\mathbf N}}

\newcommand{\GcS}{G\times\cS}
\newcommand{\GfS}{G\times\fS}

\newcommand{\fS}{{\mathfrak S}}

\newcommand{\E}[1][]{{\mathcal E^{#1}}} 
\newcommand{\ES}[1][]{{\E^{\cS}_{#1}}}
 
\newcommand{\DD}[1][]{{\mathcal D^{#1}}}

\newcommand{\bfJ}{\mathbf J}
\newcommand{\bfI}{\mathbf I}

\begin{document}   

\title 
{Quantization of  Double Covers of Nilpotent Coadjoint Orbits I: 
Noncommutative Models} 
\author{Ranee Brylinski}
\address{Department of Mathematics,
        Penn State University, University Park 16802}
\email{rkb@math.psu.edu}
\urladdr{www.math.psu.edu/rkb}
\thanks{Research supported in part  by NSF  Grant No. 
DMS-9505055}
 
%\date{}

\begin{abstract}
We construct by geometric methods a noncommutative model 
$\E$ of the algebra   of regular functions on the 
universal ($2$-fold) cover  $\M$ of  certain  nilpotent coadjoint orbits  
$\OO$ for a complex simple Lie algebra $\g$. Here
$\OO$ is the dense orbit in the cotangent bundle of the generalized flag 
variety $X$ associated to a complexified Cartan decomposition
$\g=\pn\oplus\kk\oplus\pp$ where $\ppn$ are Jordan algebras by the   
TKK construction.  

We obtain $\E$ as the algebra of $\g$-finite differential
operators on a smooth Lagrangian subvariety in $\M$  
 where $\g$  is given by differential operators $\pilao^x$   twisted
according to a critical parameter $\lao=\half\pm\frac{1}{4m}$.
After  Fourier transform,  $\E$ is  a ``quadratic" extension of the algebra
$\Dlao(X)$ of  twisted differential operators for the (formal) $\lao$th power  
of the canonical bundle.

Not only is $\E$  a Dixmier algebra for $\M$, in the sense of the 
orbit method, but also
$\E$ has a   lot of additional structure, including an anti-automorphism,
a supertrace, and a non-degenerate supersymmetric bilinear pairing.
We show that $\E$ is
the specialization at $t=1$ of a graded (non-local) 
equivariant star product with parity. 
\end{abstract}
 
\maketitle      
 
{\footnotesize
{\bf Contents}   

\S\ref{sec1}. Introduction 

\S\ref{sec2}.  {The orbit cover   $\M$}  

\S\ref{sec3}. {Noncommutative models of   $\RO$} 
  
\S\ref{sec4}. {Building a  noncommutative model  of $\RM$}  

\S\ref{secE}. {The algebras $\Ela$ and symplectic geometry of $\M$}  

\S\ref{sec6}.  {Critical values of  $\la$}  

\S\ref{sec7}. {The noncommutative model  $\Elao$}

\S\ref{sec8}. {Dixmier product on $\RM$}}   
  
\section{Introduction} 
\label{sec1}      
Let $G$ be a simply-connected  semisimple Lie group with
maximal   compact subgroup $U$.  Let $\OO$ be a coadjoint orbit of  $G$
which is  stable under  dilations; so $\OO$ is
\emph{nilpotent}, i.e., identifies with an adjoint orbit of nilpotent
elements. Let $\M$ be a  Galois cover of $\OO$ with Galois
group $\cS$. Then $\M$ is a real symplectic manifold with KKS
symplectic form $\om$  and $\M$ is a  Hamiltonian $G$-space.
According to the orbit method and geometric quantization,
quantization of $\M$ should give  a unitary representation $E$ of
$G\times\cS$ whose  $\cS$-isotypic components are irreducible for $G$. 

Now suppose $G$ is in fact a complex group. Then $\M$ has 
additional  geometric structure: $\M$ is a complex symplectic manifold
with  complex structure $\bfI$ and holomorphic KKS form
$\Om$ where $\mathrm{Re}\,\Om=\om$. Moreover the data  $(\bfI,\Om)$
extend  to a $U$-invariant hyperkahler structure (\cite{Kr}). We 
now expect additional  structure for the quantization.
First  $E$ should be realized as a Hilbert  space  of $\bfI$-holomorphic
functions on  $\M$. (Indeed, the hyperkahler data provides a
``new" complex structure $\bfJ$   which polarizes $\om$. Then
$E$ should consist   of $\bfJ$-holomorphic functions. Finally we can 
``rotate"   $\M$ to move $\bfJ$ to $\bfI$.)
Second the Harish-Chandra module 
$\E$ of $U$-finite vectors in $E$ should be, as a representation of $\GcS$,
the  algebra of $U$-finite $\bfI$-holomorphic functions on  $\M$; notice  this
is  the algebra   $\RM$ of $\bfI$-holomorphic functions  which
are  regular in the  sense of algebraic geometry.

Third, and most importantly, $\E$ should be  a
$(\GcS)$-equivariant  \emph{noncommutative model} of  $\RM$.
This means that $\GcS$ acts on $\E$ by algebra automorphisms and 
$\E$ has an invariant algebra filtration for which
$\gr\E$  is  equivariantly isomorphic, as a graded Poisson algebra, to
$\RM$.   Here $\RM$ has the Euler grading extending over $\halfN$
(see  \cite{Mo}, \cite{B-K:jams}). 
In particular, $\ES$ is a noncommutative model of $\RO$.  
Let $\g$ be the Lie algebra of $G$. The
Hamiltonian functions $\phi^x\in\RO$, $x\in\g$, 
(defined by the embedding of $\OO$ into $\g^*$)
lift to elements $\psi^x\in\ES$  such that $[\psi^x,\psi^y]=\psi^{[x,y]}$.
Then   $\ES$ (or at least  a subalgebra)  is   generated by the $\psi^x$  and
is a  primitive  quotient $\Ug/J$ of the  universal enveloping algebra.

This additional structure suggests that  a reasonable first objective of
quantization is to construct an $\cS$-equivariant noncommutative model
$\E$ once we are given  $\M$.  Then the second objective is to  establish
unitarity using  the model; see   Remark \ref{rem:pos} 
and   \cite{me:2kpos}. The case where the cover $\M$ is non-trivial, i.e.,
$\M\neq\OO$, is particularly intriguing as then $\E$ must be  
strictly larger  than   $\ES$ and so strictly larger  than  $\Ug/J$. 

A noncommutative model $\E$ is an example of a \emph{Dixmier algebra},
that is   an overring, with certain finiteness properties,  of a  primitive
quotient of $\Ug$.     Vogan (\cite{Vog1}-\cite{Vog4}) has
developed  this approach to quantization of $\M$  (at least the first
objective)  as a natural extension of  the   orbit method in the Dixmier 
sense (a correspondence between primitive ideals  and orbits).   
Joseph has developed Dixmier algebra theory from a more ring-theoretic
perspective; see \cite{Jos} and references therein, especially 
the paper \cite{J-S} with Stafford. For results particularly related to our
paper, see also \cite{McG1}-\cite{McG4}, \cite{Mo}, \cite{Za}.    

In this paper, we construct by geometric methods a noncommutative
model $\E$ for the universal cover $\M$ in a family of cases
where  $\M$   has degree $2$. This
quantization has some very nice properties and demonstrates, in a rather
non-trivial way,   some basic  paradigms of   quantization.  See in
particular  Corollary \ref{cor:3} and Theorem \ref{thm:4}.
Not only is $\E$  a Dixmier algebra for $\M$  but also $\E$ has a  
lot of additional structure, including an anti-automorphism,
a supertrace, and a non-degenerate supersymmetric bilinear pairing. 

The nilpotent orbits we work arise in the following way
(see \S\ref{sec2} and especially Table \ref{table:O}). We assume
$\g$ is simple and has a symmetric  subalgebra $\kk$ with   one-dimensional
center so that $\g=\pn\oplus\kk\oplus\pp$.
We further assume that there is a non-constant homogeneous primitive
$\kk$-semi-invariant polynomial function $F$ on $\pn$;  then $\pn$ is a 
Jordan  algebra  by the TKK construction
and $F$ is its Jordan norm.  We take  $\OO$ to
be the orbit of a generic element of $\pn$ and then (on account of $F$) the
universal cover $\M$ is $2$-fold and so $\cS=\Z_2$. 

In fact  $\OO$ has a more geometric realization as the
unique Zariski open dense $G$-orbit in the cotangent bundle $T^*X$ of a
(generalized) flag variety $X$. Here $X$ is $G/Q^-$ where
$Q^\pm$ is the parabolic subgroup with Lie algebra
$\q^\pm=\kk\oplus\ppn$. See \S\ref{ss2_mcO} and, for an example, 
\S\ref{ss2_surr}. A nice fact is that all regular functions on    
$\OO$ extend to  $T^*X$.

Each algebra  $\RO=\RTX$ has a natural family of noncommutative
models which we study in \S\ref{sec3}. These are the algebras
$\Dla(X)=\D(X,\NN^\la)$ of  twisted differential  operators,
equipped with the order filtration,  where $\NN$ is
the canonical line bundle and $\la$ is any complex number.  (If
$\la$ is not integral then  $\NN^\la$ must be interpreted formally.)
The   functions  $\phi^x$ lift to   twisted vector fields $\etala^x\in\Dla(X)$ 
which  generate $\Dla(X)$  so that $\Dla(X)=\Ug/\Jla$. 
For general $\la$, $\Dla(X)$ is a simple ring (see Proposition
\ref{prop:Dlax=simple})  and so $\Jla$ is a maximal $2$-sided ideal.

With this as our starting point, the problem is to build a filtered overring
of $\Dla(X)$ which is then a noncommutative model of $\RM$. 
To do this we introduce a complex algebraic Lagrangian submanifold
$\tZ$ in $\M$. Here $\tZ$ covers  (the orbital variety) 
$Z=(F\neq 0)=\OO\cap\pn$ and $\tZ$ admits the   function $w=\sqrt{F}$.

We embed $\Dla(X)$ into the algebra $\DtZ$ of 
differential operators  on  $\tZ$  by the sequence 
\begin{equation}\label{eq:seq_into} 
\Dla(X)\hookrightarrow\Dpp\xrightarrow{\FT}\Dpn\hookrightarrow\DZ
\hookrightarrow\DtZ
\end{equation} 
Here the first map is restriction to the big cell identified with $\pp$, 
the second is the Fourier transform (following \cite{Gon}),  
the third is restriction to (the Zariski
open dense set) $Z$, and the fourth map is defined by lifting differential 
operators. See  \S\ref{sec4}.  

For each value of $\la$ we   realize $\g$ inside $\Dpn$  
by the operators  $\pila^x=\FT(-\etala^x)$. 
Then $\pila^x$ is multiplication by $x$ 
if $x\in\pp$, $\pila^x$ is a twisted vector field if $x\in\kk$, or $\pila^x$
has order $2$ if $x\in\pn$. These $\pila^x$ are familiar in representation
theory because they make $\Spp$ into a lowest weight
representation, in fact a generalized Verma module for $\q^-$.
See \S\ref{ss4_pilax}.
Our first main result (Theorem   \ref{thm:Ela}) is
\begin{thm}\label{thm:1}  
Let $\Ela$ be the $\g$-finite part of $\DtZ$ 
with respect to the operators  $[\pila^x,\cdot]$.
Then $\Ela$ has a natural $\g$-stable algebra filtration
and we have canonical inclusions \break
$\RO\subseteq\gr\Ela\subseteq\RM$.
\end{thm}
We prove this in \S\ref{ssE_proof} by constructing a new filtration 
on $\DtZ$ which extends the Fourier transform of the order filtration on
$\Dpp$. We show  its symbol calculus  produces
a symplectic open embedding of $\TtZ$ into $\M$ and then we use the fact
that  $\RM$ is the $\g$-finite part of its fraction field.
 
The condition now for $\Ela$ to be a noncommutative model of $\RM$ is
that  $\gr\Ela=\RM$.   In our second main result (Theorem \ref{thm:la})
we figure out which  values of $\la$  satisfy this. This involves the
positive constant $m$ attached to  $\OO$  by the 
property that the $Q^+$-semi-invariant section in  
$\Gamma(X,\NN\inv)$ has weight $\chi^{2m}$ where $\chi$ is the weight of
$w$.    
\begin{thm}\label{thm:2} 
We have $\gr\Ela=\RM$ if and only if $\la=\half\pm\frac{1}{4m}$.
\end{thm}
To prove this, we reduce (in Corollary \ref{cor:Ela}) to showing that the
function $w$, regarded as a multiplication operator, is $\g$-finite.
Then in, \S\ref{ss6_S}-\S\ref{ss6_com}, we  use Jordan algebra techniques 
to compute how $w$ transforms under  $\U(\pn)$. We expect these Jordan
techniques belie some deeper connection. Also it would be interesting to
extend our result  to describe the full set of rings $\gr\Ela$ that appear
as we vary $\la$.

Let  $\lao=\half-\frac{1}{4m}$ and $\laop=\half+\frac{1}{4m}$.
If $\g=\fsl(2,\C)$, then     $\lao=\four$,
$\laop=\frac{3}{4}$ and  $\Elao$ is the Weyl algebra $\C[\pdb{w},w]$.  
(See  \S\ref{ss4_sl2}, Example \ref{ex:sl2_sharp} and Remark
\ref{rem:3/4}. Example \ref{ex:sl2_bqcirc}  and Remark \ref{rem:Lax}.)
We are finding a new reason for the old result attaching  $\C[\pdb{w},w]$
to $\M$.  

The two algebras $\Elao$ and $\Elaop$ are both anti-isomorphic 
and isomorphic (Corollary \ref{cor:tdelta} and Proposition \ref{prop:Innw}).
It follows (Corollary \ref{cor:Jlao})
that $\Jlao$ coincides with $\Jlaop$ and is  stable under
the principal anti-automorphism $\tau$ of $\Ug$.
Thus our construction of $\Elao$ attaches to $\M$ just one  $2$-sided ideal 
of $\Ug$, namely the maximal ideal $\Jlao$.
It is reasonable to attach to $\OO$ the maximal ideal $J^{\half}$.
We may then think of  the values $\lao$ and $\laop$  as representing 
some ``quantum fluctuation" about   $ \half$, caused by
passing from $\OO$ to $\M$.

In \cite[Tables 5-10]{McG4}, McGovern attached, by a completely
different method, a Dixmier algebra $\DD$ to
$\M$  in each of our cases (except when $\gR=\fso(2,p)$ and $p$ is odd).
McGovern starts by manufacturing an infinitesimal character  according to
the recipes formulated in his Yale preprints and \cite{McG3}. Then 
he invokes Moeglin's construction in \cite{Mo}  so that $\DD$
is the  $\g$-finite part of the endomorphism ring of a certain
(degenerate) simple Whittaker module.
McGovern conjectures that $\gr\DD$ is isomorphic to $\RM$ and says
that he  can check this case-by-case when $\g$ is classical. 
By comparing infinitesimal characters  
(Corollary \ref{cor:infch} and Remark \ref{rem:McG}) and applying 
results of Moeglin on $\gr\DD$, we find independently
that McGovern's algebra $\DD$ always coincides with our algebra $\Elao$.
In particular, this proves McGovern's conjecture in our one  case 
$\gR=\fe_{7(-25)}$  where $\g$ is exceptional.
It would be extremely valuable to find a geometric construction   of the
Whittaker module; we conjecture that this can done in the context of our 
construction of $\Elao$.

Now returning to Theorem \ref{thm:1}, we have a built-in module for
$\Elao$, namely the ring    of regular functions on $\tZ$.  
Consider the submodule
$\cH$ generated by the constant function $1$. We  prove  (Proposition
\ref{prop:cH=fs} and Corollary \ref{cor:Elao==})   
\begin{cor}\label{cor:3} 
$\cH$ is a faithful simple module for $\Elao$.  Moreover $\cH$ identifies
with a subalgebra of $\RM$ which is maximal Poisson abelian.
Then $\Elao$ is   the  $\g$-finite part, 
with respect to the operators  $[\pilao^x,\cdot]$,
of  the algebra of differential operators on $\cH$. 
\end{cor} 
 
As a $\g$-representation, $\cH$ is 
the direct sum of   two lowest weight representations 
$\Spp$ and $w\Spp$. From this point of view,  Corollary 
\ref{cor:3} is telling us how to locate  $\Elao$ inside 
the  (computable) algebra $\End\sgfin(\cH)$ by using the algebra structure
on $\cH$.
 
Our   noncommutative model $\Elao$ is  naturally $\cS$-equivariant. 
This figures into the algebraic structure of $\Elao$.
In particular the maximality of $\Jlao$ 
``induces upward" so that $\Elao$ is
a simple ring (Corollaries  \ref{cor:Ela=simple} and \ref{cor:Ela0=simple}).
See also  Corollaries \ref{cor:ann} and \ref{cor:bi_simple} for the
decomposition  of $\Elao$ as a $\Ug$-bimodule.

In \S\ref{sec8}, we focus on interpreting $\Elao$ as a quantization of
$\RM$. We find in  \S\ref{ss8_circ} a natural quantization map   
$\bq:\RM\to\Elao$.  Via $\bq$, multiplication on $\Elao$  defines a
new $(\GcS)$-invariant  product  $\circ$ on $\RM$. 
If $\phi$ and  and $\psi$ are   homogeneous of degrees $j$ and $k$,
then  $\phi\circ\psi=\sum_{p\in\bbN}C_p(\phi,\psi)$ where
$C_p(\phi,\psi)$ is homogeneous of degree $j+k-p$.
In fact, $\circ$ deforms the  Poisson algebra
structure   in the sense that
$C_0(\phi,\psi)=\phi\psi$ and  
$C_1(\phi,\psi)-C_1(\psi,\phi)=\{\phi,\psi\}$.  

Let $\T:\RM\to\C$ be  the projection to the constant term defined by the
Euler grading.  We make $\RM$ into a supervector space where $\RM[j]$
is even or odd according to whether $j\in\bbN$ or $j\in\bbN+\half$. 
Our  third main result  
(Theorem \ref{thm:circ}, Corollaries \ref{cor:star} and \ref{cor:Lax}) is  
\begin{thm}\label{thm:4} 
With respect to  $\circ$, $\RM$ is a noncommutative superalgebra with
supertrace
$\T$. The pairing $\Q(\phi,\psi)=\T(\phi\circ\psi)$ on $\RM$ is 
$(\GcS)$-invariant, supersymmetric, non-degenerate and orthogonal for
the Euler grading.

Our product $\circ$ on $\RM$ is the  specialization at 
$t=1$ of a graded strongly $\g$-invariant \textup{(}non-local\textup{)}
star product. This   has parity, i.e., 
$C_p(\phi,\psi)=(-1)^pC_p(\psi,\phi)$.  If $x\in\g$ then
$\phi^x\circ\psi=\phi^x\psi+\half\{\phi^x,\psi\}+\La^x(\psi)$ where
 $\La^x$ is  the $\Q$-adjoint of   $\psi\mapsto\phi^x\psi$.
\end{thm}
The parity condition, which is essential for star products, is not
automatic but  comes from an anti-automorphism $\be$ of 
$\Elao$ which extends the principal anti-automorphism $\tau$ of
$\Ug/\Jlao$. We find that $\be$  falls out of our comparison 
of  $\Elao$ with  $\Elaop$  (Corollary  \ref{cor:be}). 

Our star product is \emph{non-local} in the sense that the operators
$C_p(\cdot,\cdot)$  fail  in general to be bi-differential. We know this
because   already the operators $\La^x$ fail  in general to be   differential
general (see \cite{me:2kpos}  and Remark \ref{rem:Lax}).

We say $\circ$ is a \emph{Dixmier product} because it makes
$\RM$ into a Dixmier algebra for $\M$, equipped with extra structure.
In \cite{me:DixtoStar}
we show, for  a general nilpotent orbit cover $\M$,   how adding
some   axioms  (for $\be$, $\T$ and $\Q$) to the usual
definition of Dixmier algebra  produces a Dixmier product on $\RM$
and all the results in \S\ref{ss8_circ}-\S\ref{ss8_Lax}.  
Thus these results persist even when $\RM$ has non-trivial 
multiplicities.

In \cite{me:2kpos} we show that our star product on $\RM$ is ``positive"  
in a   sense which we define and consequently $\RM$ becomes a unitary
representation of $G$ (see Remark \ref{rem:pos}) made up of two
irreducible components.  

I  thank Jean-Luc Brylinski, Michel Duflo, Tony Joseph,
Siddhartha Sahi, Eric Sommers, Toby Stafford,  and David Vogan   for
helpful conversations regarding their own work, related work of others,
and the philosophy of what  quantizations of orbit covers should look
like.  I also thank  Aravind Asok, Alex Astashkevich and Francois Ziegler
for  what I learned in collaboration with them on related problems  over 
the past four  years.

\section{The orbit cover    $\M$}  
\label{sec2}   

\subsection{Momentum construction of \boldmath$\OO$} 
\label{ss2_mcO}
Let $G$ be a connected and simply-connected complex semisimple
Lie group with Lie algebra $\g$. Let $\GR$ be a real form  of $G$.
Let $K\subset G$ be the complexification of a maximal compact
subgroup $\KR$ of $\GR$.  Let $\gR$, $\g$, $\kR$, $\kk$ be
the Lie algebras of $\GR$, $G$, $\KR$, $K$. Let $\g\to\g$,
$x\mapsto\ovl{x}$, be  the complex conjugation map. Then we have the
Cartan decomposition $\gR=\kR\oplus\pR$ and its complexification
$\g=\kk\oplus\p$.

We assume from now on that  the real symmetric pair $(\gR,\kR)$, 
or equivalently the  complex symmetric pair $(\g,\kk)$, is
\emph{Hermitian}. This means that  $[\p,\p]=\kk$ and there exists
$x_0\in\Cent\kR$ such that  $\adj\,x_0$ defines a complex
structure  on $\pR$. Then we get the splitting
$\p=\pp\oplus\pn$  where   $\ppn$ are the $\pm i$-eigenspaces  
of $\adj\,x_0$, and so we  get
\begin{equation}\label{eq:g=pnkkpp} 
\g=\pp\oplus\kk\oplus\pn
\end{equation} 
Every Hermitian symmetric pair   is a direct sum, in the
obvious way, of  Hermitian symmetric pairs
$(\g,\kk)$ where $\g$ is simple; such pairs are called
\emph{irreducible}. 

Now $\ppn$   are  (complex conjugate) abelian Lie 
subalgebras of $\g$.  Let $U^{\pm}\subset G$ be the
corresponding abelian  subgroups. Then 
$Q^{\pm}=KU^{\pm}$ are parabolic subgroups of $G$ with Lie algebras 
$\q^{\pm}=\kk\oplus\p^{\pm}$.
The coset spaces $G/Q^{\pm}$ are then   (generalized) flag varieties of
$G$. We put   $X=G/Q^-$.
 
We differentiate the $G$-action on $X$ to obtain an infinitesimal action 
\begin{equation}\label{eq:etax} 
\g\lto\Vect\,(X),\qquad x\mapsto\eta^x
\end{equation}
where the value of $\eta^x$ at a point $q$ is
$\eta^x_q=\frac{\dd\phantom{t}}{\dd{t}}|_{t=0}(\exp -tx)\cdot q$. 
Here $\Vect(X)$ denotes the Lie algebra of algebraic holomorphic   
vector fields on $X$ and by \emph{infinitesimal action} we mean that
(\ref{eq:etax})   is a Lie algebra homomorphism.
 
The natural action of $G$ on $X$ induces a canonical  action
of $G$ on the  cotangent bundle $T^*X$. This $G$-action is Hamiltonian
with respect to the  canonical  symplectic form 
$\Om$ on $T^*X$; throughout this paper  \emph{symplectic} means 
\emph{algebraic holomorphic symplectic}.
The $G$-equivariant  moment map
\begin{equation}\label{eq:mu} 
\mu:T^*X\lto\g^*
\end{equation}
is defined by   $\langle\mu(m), x\rangle=\mu^x(m)$ where
$\mu^x$ is the order one symbol of $\eta^x$. In other words, 
the comorphism  
\begin{equation}\label{eq:mux} 
\mu^*:\g\lto\RTX,\qquad x\mapsto\mu^x
\end{equation}
is a  Lie algebra homomorphism where  $\RTX$ is 
equipped with the Poisson bracket $\{\cdott\}$ defined by $\Om$.  
Thus $\mu^*$ defines  Hamiltonian  $\g$-symmetry  on $T^*X$.
If we identify $T^*X$ with the contracted product bundle  
$G\times_{Q^-}(\g/\q^-)^*$ in the usual way, then $\mu$ is the
collapsing map. The following   fact defines $\OO$ for us.
\begin{prop}\label{prop:OO} 
The image of the moment map
$\mu$ is the closure of a single nilpotent coadjoint orbit 
$\OO$ in  $\g^*$ so that
\begin{equation}   
\Cl(\OO)=\mu(T^*X)=G\cdot(\g/\q^-)^*
\end{equation}
Then $\OO$ is the Richardson orbit  associated to  $Q^-$. 
The map $\mu$ is generically $1$-to-$1$ and $\Cl(\OO)$ is normal.
Consequently all   regular 
functions on $\OO$  extend to $\Cl(\OO)$ so that  $\RO=R(\Cl(\OO))$.

The moment map $\mu$ is bijective over $\OO$ and $\mu\inv$ defines 
a  $G$-equivariant Zariski open embedding of complex algebraic
manifolds
\begin{equation}\label{eq:bj} 
\bj:\OO\lto T^*X 
\end{equation}
Let $\om$ be the  KKS symplectic form  on  $\OO$. 
Then $\bj$ is   symplectic; i.e., $\bj^*\Om=\om$.  
The comorphism  $\bj^*:\RTX\lto\RO$
is an isomorphism of Poisson algebras. 
\end{prop}

To say  the coadjoint orbit $\OO$ is \emph{nilpotent}
means that if we identify $\g^*$ with $\g$ using the complex 
Killing form $\Killg{\cdott}$ of $\g$, then the adjoint orbit 
corresponding  to $\OO$ consists of nilpotent elements. 
This happens if  and only if
$\OO$ is  stable under  the    dilation action of $\C^*$ on $\g^*$.
We often identify $\OO$ with its corresponding adjoint orbit. 

Let $\phi:\OO\to\g^*$ be the inclusion with comorphism 
\begin{equation}\label{eq:phix} 
\phi^*:\g\to \RO,\qquad x\mapsto \phi^x
\end{equation}
Then $\{\phi^x,\phi^y\}=\phi^{[x,y]}$ where the Poisson bracket
on $\RO$ is defined by the KKS form $\om$; this property
determines $\om$   uniquely.  
So $\phi^*$ is  a  Lie algebra homomorphism and thus defines 
Hamiltonian   $\g$-symmetry on $\OO$. 
 
\subsection{The graded Poisson algebra  \boldmath $\RO$} 
\label{ss2_RO}
We will say a (complex) commutative algebra $\A$   is 
\emph{graded} if $\A$ is equipped with a vector space grading
$\A=\oplus_{p\in\halfN}\,\A^p$
such that $\A^p\A^q\subseteq\A^{p+q}$.
Here $\bbN=\{0,1,2,\dots\}$.
In practice,  $\A$ will be a finitely generated algebra with $\A^0=\C$.
 
\begin{definition}\label{def:gPa} 
A \emph{graded Poisson algebra} is a Poisson algebra $\A$ together
with an algebra grading $\A=\oplus_{p\in\halfN}\,\A^p$ such that
$\{\A^p,\A^q\}\subseteq \A^{p+q-1}$. If also $\A^p=0$ when
$p\notin\bbN$, then we say $\A$ is an 
\emph{$\bbN$-graded} Poisson  algebra.
\end{definition}
We have three examples  on hand of $\bbN$-graded Poisson algebras.
(i) The symmetric  algebra  $\Sg$   with  grading defined by polynomial
degree and   Poisson bracket defined by the Lie bracket on $\g$.
(ii) $\RTX$ with  Poisson bracket defined by $\Om$ and   grading  
$\RTX=\oplus_{p\in\bbN}\,\RTX[p]$ where
$\RTX[p]$ is the subspace of functions which are
homogeneous of degree $p$  on the fibers of  $T^*X$ over $X$.
(iii) $\RO$ with Poisson bracket defined by $\om$ and  Euler grading
$\RO=\oplus_{p\in\bbN}\,\RO[p]$ where $\RO[p]$ is the subspace of 
homogeneous degree $p$ functions.

The natural extensions $\mu^*:\Sg\to\RTX$, $P\mapsto\mu^P$, and 
$\phi^*:\Sg\to\RO$, $P\mapsto\phi^P$,
are graded  Poisson  algebra homomorphisms. Plainly   $\mu\bj=\phi$
and so Proposition \ref{prop:OO} gives
\begin{cor}\label{cor:I} 
The two maps $\mu^*$ and $\phi^*$ are surjective with the same
kernel   $I\subset S(\g)$. Then $I$ is 
the ideal of functions vanishing on $\OO$.  
\end{cor}
 
\subsection{The tube condition} 
\label{ss2_tube}

We define $h\in\kk$ by   $h=-ix_0$
so that $\ppn$ is the $\pm 1$-eigenspace  of $\adj h$.
A  Hermitian symmetric pair
$(\g,\kk)$ is said to be of \emph{tube type} if there exists
$e\in\pp$ such that    $\s=\C e+\C h+\C\ove$
is a Lie subalgebra of $\g$ isomorphic to $\fsl(2,\C)$ with
bracket relations: $[h,e]=e$, $[h,\ove]=-\ove$, $[e,\ove]=2h$. 
In this case, $e$ lies in $\OO$ and so $\OO=G\cdot e=G\cdot\ove$.
Clearly $(\g,\kk)$ is of tube type if and only if every one of its 
irreducible components is of tube type. 

\begin{lem}\label{lem:pi1O} 
Assume that $\g$ is simple. Then the fundamental group  
of $\OO$  is  $\Z_2$ if $(\g,\kk)$ is of tube type or is trivial otherwise.
\end{lem}

In Table \ref{table:O}, we write down the familiar list of  irreducible
Hermitian symmetric pairs $(\g,\kk)$ of tube type
along with the real form $\gR$, the  rank
$r$ of $(\g,\kk)$, and $n=\dim\ppn$. Note that
$\dim\OO=2\dim X=2n$.
We specify the orbit $\OO$. 
In every case except one, the subalgebra $\kk\subset\g$ is unique up
to conjugacy and hence we get a single orbit $\OO$. The exception is
the case where $\gR=\fso^*(4r)$, as then there are two choices for
$\kk$ and these give rise to two distinct orbits which are  exchanged
by  outer automorphism. For the classical cases we give the
partition indexing $\OO$ (see e.g. \cite{C-M}) and in the exceptional
case we give the dimension of $\OO$ (this is enough since there is only
one nilpotent orbit of that dimension). 
In Table \ref{table:O}, $r\ge 1$ and $p\ge 2$.

\begin{table}[h] 
\caption{The Orbits $\OO$}
\label{table:O}
\begin{equation}\nonumber 
\begin{tabular}{|c|c|c|c|c|c|}  
\hline \rule[-3mm]{0mm}{8mm}
$\g$&$\kk$&$r$&$n$&$\gR$&$\OO$\\ 
\hline \rule[-3mm]{0mm}{8mm}
$\fsp (2r,\C)$&$\fgl(r,\C)$&$r$&$\half r(r+1)$&
$\fsp(r,\R)$&$(2^r)$\\[3pt]
$\fsl(2r,\C)$&$\s(\fgl(r,\C)\oplus\fgl(r,\C))$&$r$&$r^2$&
$\fsu(r,r)$&$(2^r)$\\[3pt]
$\fso(4r,\C)$&$\fgl(2r,\C)$&$r$&$r(2r-1)$&$\fso^*(4r)$&
$(2^{2r})_{I,II}$\\[3pt]
$\fe_7$&$\fe_6\oplus\C$&$3$&$27$&
$\fe_{7(-25)}$&$54$\\[3pt]
$\fso(2+p,\C)$&$\fso(2,\C)\oplus\fso(p,\C)$&$2$&$p$&$\fso(2,p)$
&$(3,1^{p-1})$\\[3pt]
\hline
\end{tabular} 
\end{equation}
\end{table} 

The geometric interpretation of the tube condition is that the
Hermitian symmetric space $\GR/\KR$ is of tube type. 
The  Jordan theoretic interpretation is that in TKK theory, $\ppn$ is
not just a  Jordan triple system but also a Jordan algebra; cf.  
\S\ref{ss6_Jordan}. The invariant theoretic  interpretation is given
by the next lemma in terms of the algebra $\Spp^{K-semi}$ of
$K$-semi-invariants.
\begin{lem}\label{lem:F} 
Assume $\g$ is simple. Then   $\Spp^{K-semi}\neq\C$
if and only if $(\g,\kk)$ is of tube type. In the tube case,
$\Spp^{K-semi}=\C[F]$ is a polynomial ring 
in one  homogeneous generator $F$. Then   $F$ has degree  $r$ 
and $F$ is unique up to scaling. The weight of $F$  is
$\chi^2$ where $\chi$ is a generator of the character group  of $K$. 
\end{lem}
This lemma defines $\chi$ and then we extend $\chi$ to a character  
of $Q^\pm$ which is trivial on the unipotent radical. 

\subsection{The universal cover \boldmath $\M$} 
\label{ss2_tO}
From now on we assume that $(\g,\kk)$  is an irreducible Hermitian
symmetric pair  of tube type.  Then, by Lemma \ref{lem:pi1O}, $\OO$
admits a universal  $2$-fold   covering $\M$. 
We can give a nice geometric construction of $\M$ using 
Lemma \ref{lem:F}  and the homogeneous function  $\phi^F\in\RO[r]$ 
defined by the function $F$ introduced in Lemma \ref{lem:F}.
\begin{prop}\label{prop:tO} 
The function $\phiF$ is not a square in  the field $L=\C(\OO)$ of
rational functions on $\OO$. Let $\tL$ be the field extension of $L$
defined   by adjoining
\begin{equation}\label{eq:ze=} 
\ze=\sqrt{\phiF}
\end{equation} 
The $G$-representation  on $L$ extends uniquely to $\tL$ and 
$\ze$ is $Q^+$-semi-invariant of weight $\chi$.  

The normalization of $\OO$ in $\tL$ is a $G$-homogeneous complex
algebraic manifold $\M$. The normalization map is a 
$G$-equivariant $2$-fold covering $\kap:\M\to\OO$. 
$\RM$ is the algebra of $G$-finite functions in $\tL$. Each
irreducible $G$-representation occurring in $\RM$ has multiplicity one 
and is self-dual.  

The function $\ze$  lies in $\RM$ and is the highest weight
vector   of  a   finite-dimensional irreducible $G$-representation 
$V\subset\RM$; $V$ is given in Table \textup{\ref{table:V}}.
The algebra $\Z_2$-grading  defined by the action of the Galois group  
$\cS=\Z_2$ is 
\begin{equation}\label{eq:RtO=} 
\RM=\RO\oplus\RO V
\end{equation}
\end{prop}

\emph{Explanation of Table \textup{\ref{table:V}}:}
In the first row, $\wedge_o^r(\C^{2r})$ is the kernel of the map
$\wedge^r(\C^{2r})\to \wedge^{r+2}(\C^{2r})$ defined by taking
the wedge product with the symplectic form. In the third row, we
obtain the the two half-spin representations, corresponding to the
two orbits listed in Table \ref{table:O}. The other entries are clear.

\begin{table}[h] 
\caption{The  Representation $V$}
\label{table:V}
\vskip 12pt
\center{
\begin{tabular}{|c|c|c|}
\hline \rule[-3mm]{0mm}{8mm}
$\gR$&$V$\\ \hline  \rule[-3mm]{0mm}{8mm}
$\fsp(r,\R)$&$\wedge_o^r(\C^{2r})$\\[2pt]
$\fsu(r,r)$&$\wedge^r(\C^{2r})$\\[2pt]
$\fso^*(4r)$&$\C^{2^{2r-1}}=\pm\half\mbox{-spin}$\\[2pt]
$\fe_{7(-25)}$&$\C^{56}$\\[2pt] 
$\fso(2,p)$&$\C^{2+p}$\\[2pt]
\hline 
\end{tabular}} 
\end{table}

We have several easy consequences of Proposition \ref{prop:tO}.
\begin{cor}\label{cor:tO} 
The orbit $G\cdot(e,\ze)$ inside  $\OO\times V$  is a
$G$-equivariant model for $\M$ where the composition 
$G\cdot(e,\ze)\hhookrightarrow\OO\times V\llongrightarrow\OO$ 
gives the $2$-fold covering onto $\OO$.
\end{cor}
\begin{remark}\label{rem:tO_sl2}  
In just  one case, namely $\g=\fsl(2,\C)$, the  map
$\M\hookrightarrow\OO\times V\to V$ is  an embedding.
Here $V=\C^2$. This accounts for why the $\fsl(2,\C)$ is so easy to write
down; see \S\ref{sec1} and \S\ref{ss2_surr}.
\end{remark}

\begin{cor}\label{cor:RtO=PA} 
The KKS form $\om$ lifts to a $G$-invariant  symplectic form 
$\tom=\kap^*\om$ on $\M$ which then defines a Poisson bracket
on $R(\M)$ and  $\C(\M)$. 
The $G$-representation   on $\RM$ corresponds to the $\g$-representation
\begin{equation}\label{eq:Phi}  
\Phi:\g\to\End\,\RM,\qquad  \Phi^x=\{\phi^x,\cdot\}  
\end{equation}
As a Poisson algebra, $\RM$  is generated by $\RO$ and   $\ze$. 
\end{cor}

\begin{cor}\label{cor:cG} 
Let $\cG\subset\RM$ be the set of functions which
Poisson commute with $\phi^x$ for all $x\in\pp$.
Then $\cG$ is a maximal Poisson abelian  subalgebra of $\RM$. 
We have $\cG=\{\phi^T+\phi^{T'}\ze\,|\, T,T'\in\Spp\}$.  
\end{cor}
\begin{remark}\label{rem:corcG}  
The $K$-types in $\cG$ are in natural bijection with the $G$-types in $\RM$.
\end{remark}
       
The square of the $\C^*$-action on $\OO$ lifts to the $\C^*$-action on
$\M$ defined by $s(u,v)=(s^2u,s^rv)$  in the model of Corollary \ref{cor:tO}.
Notice that $-1$ interchanges points in the cover 
$\kap:\M\to\OO$ if $r$ is odd, or acts trivially if $r$ is even.  Let
$\RM[j]\subset\RM$ be the space of homogeneous degree $2j$
functions where $j\in\halfZ$.  Then $\RO[p]\subseteq\RM[p]$ for
$p\in\bbN$.   
\begin{cor}\label{cor:RtO=oplus}
$\RM$ is  a graded  Poisson algebra with respect to the  
Euler grading $\RM=\oplus_{j\in\halfN}\,\RM[j]$.  
We have $V\subseteq\RM[\frac{r}{2}]$.
\end{cor}
\subsection{Example: \boldmath$\gR=\fsu(r,r)$} 
\label{ss2_surr}  
Here $\g=\fsl(2r,\C)$ and
complex conjugation on $\g$ is the map
$\left(\Atop{A}{C}\Atop{B}{D}\right)
\mapsto \left(\Atop{-A^*}{B^*} \Atop{C^*}{-D^*}\right)$
where $A,B,C,D$ are complex $r\times r$ matrices. We can choose 
$\kR$ so that 
\begin{equation}\nonumber \textstyle
\kk=\left\{\left(\Atop{A}{0}\Atop{0}{D}\right)\,\big|\,
A,D\in\fgl(r,\C),\; \Tr(A+D)=0\right\} 
\end{equation}
We can pick $x_0$ (which is unique up to sign) so that
\begin{equation}\nonumber \textstyle
\pp=\left\{\left(\Atop{0}{0}\Atop{B}{0}\right)\,\big|\,
B\in\fgl(r,\C)\right\},\qquad
\pn=\left\{\left(\Atop{0}{C}\Atop{0}{0}\right)\,\big|\,
C\in\fgl(r,\C)\right\}   
\end{equation}

Then $X$ identifies with the Grassmannian $Gr(r,2r)$ of
$r$-dimensional vector subspaces $L$ in $\C^{2r}$.
A point in $T^*X$ corresponds to a pair  $(p,L)$ where
$p$ is a linear transformation $\C^{2r}/L\to L$. 
The moment map $\mu:T^*X\to\g^*$ is given by
$\mu(p,L)=\Tr(xy_{p,L})$ where $y_{p,L}$ is the composite map
$\C^{2r}\to\C^{2r}/L\mapright{p}L\hookrightarrow\C^{2r}$.
Thus $\OO$ identifies with the  nilpotent orbit
$\{x\in\fsl(2r,\C)\,|\, x^2=0,\;\rank x=r\}$.

We  have $2h=\left(\Atop{I}{0}\Atop{0}{-I}\right)$ 
where $I$ is the $r\times r$ identity matrix.  We can choose
$e=\left(\Atop{0}{0}\Atop{I}{0}\right)$ and then
$\ove=\left(\Atop{0}{I}\Atop{0}{0}\right)$.
The polynomial $F\in S^r(\pp)$ is the determinant and so
$\chi\left(\Atop{A}{0}\Atop{0}{B}\right)=\Det A$ and
$\phiF\left(\Atop{A}{C}\Atop{B}{D}\right)=\Det C$.  

\section{Noncommutative models of   $\RO$}   
\label{sec3}

\subsection{Noncommutative models}
\label{ss3_ncm} 
We will say a noncommutative algebra $\B$ is \emph{filtered} if   
$\B$ is equipped with an increasing filtration
$\B=\cup_{j\in\halfN}\,\B_j$
such that $\B_j\B_k\subseteq\B_{j+k}$.
The   associated graded algebra  is
$\gr\B=\oplus_{j\in\halfN}\gr_j\B$ where 
$\gr_j\B=\B_j/\B_{j-\half}$.
Let $\bp_j:\B_j\to\gr_j\B$ be the natural projection.
Suppose we have, for all $j,k,\in\halfN$,
\begin{equation}\label{eq:3.1} 
[\B_j,\B_k]\subseteq\B_{j+k-1} 
\end{equation}
Then $\gr\B$ is commutative and moreover $\gr\B$ is a graded
Poisson algebra  with Poisson  bracket   given by
$\{\bp_j(b),\bp_k(c)\}=\bp_{j+k-1}(bc-cb)$
where $b\in\B_j$ and $c\in\B_k$.

\begin{definition}\label{def:ncm} 
Let $\A=\oplus_{p\in\halfN}\,\A^p$ be a graded  Poisson
algebra as in Definition \textup{\ref{def:gPa}} with Hamiltonian
$\g$-symmetry given by  a Lie algebra embedding
$\phi:\g\to\A^1$, $x\mapsto\phi^x$. 
A \emph{noncommutative   model} of  $(\A,\phi)$ is a triple
$(\B,\ga,\psi)$ where $\B$ is a noncommutative  filtered algebra 
$\B=\cup_{j\in\halfN}\,\B_j$ satisfying  \textup{(\ref{eq:3.1})}, 
$\ga:\gr\B\lto\A$ is a graded Poisson  algebra isomorphism, and
$\psi:\g\to\B_1$, $x\mapsto\psi^x$,
is a   Lie algebra homomorphism 
such that $\ga(\bp_1(\psi^x))=\phi^x$ for all $x\in\g$. 
\end{definition}

We have representations of $\g$ on $\B$ and $\A$ given by the
operators  $b\mapsto[\psi^x,b]$ and $a\mapsto\{\phi^x,a\}$; 
the former induces a $\g$-representation  on $\gr\B$.
Clearly $\ga$ is $\g$-equivariant and we have 
\begin{lem}\label{lem:BisoA} 
Suppose $\A^j$ is finite-dimensional for each $j\in\halfN$.
Then  $\B$ is isomorphic to $\A$ as a $\g$-representation with
$\B_j\simeq\oplus_{k=0}^j\A^k$.
\end{lem}

If $\A$ and $\B$ are graded and filtered over $\bbN$, 
then we may form a
$\halfN$-grading of $\A$ and a $\halfN$-filtration of $\B$
by putting $\A^{p+\half}=0$ and $\B_{p+\half}=\B_p$ for 
$p\in\bbN$.  In this way the
$\bbN$-graded/filtered theory is subsumed in the 
$\halfN$-graded/filtered theory.

We  often speak of a noncommutative model of $\A$ where
$\phi$ is implicitly understood.
In particular, if $\A$ is $\RO$ or $\RM$, then we always take
the Hamiltonian symmetry to be the one defined in  (\ref{eq:phix}).

If we identify $\g$ with the diagonal in $\gog$, then
our $\g$-representation on $\B$ extends to the representation
\begin{equation}\label{eq:Pi} 
\Pi:\gog\to\End\B,\qquad \Pi^{(x,y)}(b)=\pi^xb-b\pi^y 
\end{equation}
This is a key aspect of noncommutative models.

\subsection{The algebras \boldmath 
$\Dla(X)$ of twisted differential operators}
\label{ss3_tdo} 
The canonical bundle $\NN$ on $X$ is the algebraic holomorphic
complex line bundle  given by the top exterior power of $T^*X$. 
Every $G$-homogeneous line  bundle  over $X$ is a 
(rational) tensor power  of $\NN$; this follows by Lemma \ref{lem:F}.

We can construct the sheaf $\Dla=\D_{X,\NN^\la}$  of
$\NN^\la$-twisted  differential operators on $X$ where $\la$ is any
complex number. This is a sheaf of noncommutative algebras.
When $\la$ is an integer,  the line bundle $\NN^\la$ exists and then
$\D_{X,\NN^\la}$ is the usual sheaf constructed using $\NN^\la$.
In particular $\D=\D^0$ is the usual sheaf of differential operators.

For the theory of twisted differential operators on  flag
varieties  (and more generally on algebraic manifolds) with applications
to representation theory,   see e.g., 
\cite{Be-Be}, \cite{Bo-Br},  \cite{Bj}, \cite{Ka}, \cite{Mi},
\cite{Vog1},  \cite{Vog4}.

We have a sheaf filtration 
$\Dla=\cup_{p\in\bbN}\,\Dla_p$ where $\Dla_p$ is the
subsheaf of differential operators of order at most $p$.
We have $[\Dla_p,\Dla_q]\subseteq\Dla_{p+q-1}$ and so
$\ugr\Dla$ is a sheaf  of graded Poisson algebras which is isomorphic
by the symbol map $\bsla:\ugr\Dla\to S(\TT)$  to the symmetric algebra
$S(\TT)$   of the tangent sheaf $\TT$ of $X$.

The symbol map defines a graded Poisson algebra inclusion
\begin{equation}\label{eq:bsU} 
\bsla_U:\gr\,\Dla(U)\hhookrightarrow\Gamma(U,S(\TT))=R(T^*U)
\end{equation}
where  $U$  is  Zariski open in $X$; we omit the subscript $U$ when the
context is clear. We have a natural Lie algebra homomorphism
\begin{equation}\label{eq:etalax} 
\etala:\g\lto\Vect\,(X)\lto\Dla_1(X),
\qquad x\mapsto\eta^x\mapsto\etala^x 
\end{equation}
where   $\etala^x$ is the Lie derivative $\cL_{\eta^x}$ acting on  
$\la$-twisted forms. We say $\etala^x$ is a 
\emph{twisted vector field on $X$}. Let
$\bvla$ denote the composite map
$\gr\,\Dla(X)\xrightarrow{\bsla}\RTX\xrightarrow{\bj^*}\RO$.
Now we know (see \cite{Vog1}): 
\begin{prop}\label{prop:DlaX=ncm} 
$(\Dla(X),\bvla,\etala)$, is a  
noncommutative model of $\RO$.
\end{prop}
\begin{proof} 
As $X$ is a generalized flag variety,  the sheaf   cohomology 
$H^1(X,S^p(\TT))$ vanishes and it follows that $\bsla_X$ is an
isomorphism -- see \cite[\S1, Lemma 1.4]{Bo-Br}. Hence 
$(\Dla(X),\bsla_X,\etala)$ is a  noncommutative model of
$\RTX$. This implies the result for $\RO$. 
\end{proof} 

We will use later the following result (true for any flag variety $X$).
\begin{prop}\textup{\cite{Br-Br}}\label{prop:Dlax=simple} 
$\Dla(X)$ is a simple ring if $\la$ satisfies 
$2\la\notin\Z-\{1\}$.
\end{prop}

\subsection{Relations with the enveloping algebra \boldmath $\Ug$}
\label{ss3_Ug} 

Now (\ref{eq:etalax}) extends to an algebra homomorphism
\begin{equation}\label{eq:etala} 
\etala:\Ug\to\Dla(X) 
\end{equation}
where $\Ug$ is the universal enveloping algebra of $\g$. We have the
standard algebra filtration $\Ug=\cup_{p\in\bbN}\Ug[p]$ and then 
$\etala$  is a filtered map, i.e., $\etala(\Ug[p])\subseteq\Dla_p(X)$. 

Let $\Jla$ be the kernel of  $\etala$ so that $\Jla$
is a two-sided ideal  in $\Ug$.
According to   \cite{Bo-Br} (which applies   since the moment map 
$\mu:T^*X\to\Cl(\OO)$ is birational with normal image),  the map
(\ref{eq:etala}) is surjective in each filtration degree.
Although only the untwisted case was treated in \cite{Bo-Br},
their method of proof (symbols), and hence their
result, extends immediately to the twisted case. Thus we find
\begin{prop}\label{prop:DlaX_gen} 
The algebra $\Dla(X)$ is generated by the twisted vector
fields $\etala^x$, $x\in\g$. Moreover, $\etala$    induces
a  filtered algebra isomorphism 
\begin{equation}\label{eq:etalas} 
\etalas:\Ug/\Jla\lto\Dla(X)
\end{equation} 
Hence  $\gr\Jla=I$ and $\gr\etalas:\Sg/I\to\RTX$ coincides with 
the isomorphism $\mus$. 
\end{prop}
\begin{cor}\label{cor:Jla} 
$\Jla$ is completely prime. Moreover if $2\la\notin\Z-\{1\}$
then  $\Jla$ is maximal.
\end{cor}
\begin{proof} 
The first statement follows as  $I$ is   prime and the second
by Proposition \ref{prop:Dlax=simple}.
\end{proof}

\subsection{The anti-symmetry \boldmath $\la\mapsto(1-\la)$}
\label{ss3_1-la} 
The following anti-symmetry will be important throughout the paper.
\begin{prop}\label{prop:theta} 
There is  a unique map
$\theta:\Dla(X)\to\D^{1-\la}(X)$  
such that $\theta$ is an algebra anti-isomorphism
and $\theta(\etala^x)=-\eta_{1-\la}^x$. 
\end{prop}
\begin{proof} 
Let $U$ be any Zariski open affine in $X$. Then $\Dla(U)$ is
generated by the  multiplication operators $f\in R(U)$ and the  
order  $1$ operators  $\cL_{\eta}$ where $\eta\in\Vect(U)$. 
We obtain an algebra anti-isomorphism
$\theta_U:\Dla(U)\to\D^{1-\la}(U)$
by assigning $\theta_U(f)=f$ and 
$\theta_U(\cL^{\la}_{\eta})=-\cL^{1-\la}_{\eta}$.
This follows  by checking the relations among  our generators 
of $\Dla(U)$;   see \cite[proof of Prop. 5.6.2]{AB:starmin}.  
In particular, $\theta_U(\etala^x)=-\eta_{1-\la}^x$. 

Now it is easy to see that the maps $\theta_U$ patch together to
define an anti-isomorphism 
$\btheta:\Dla\to\D^{1-\la}$ of sheaves of
algebras. Then $\btheta$ evaluated on global
sections gives  $\theta$ and we have
$\theta(\etala^x)=-\eta_{1-\la}^x$. Finally, $\theta$ is unique
since the vector fields $\etala^x$ generate $\Dla(X)$.
\end{proof}
 
Let $\tau:\Ug\to\Ug$ be the  algebra  anti-automorphism of
$\Ug$ such that $\tau(x)=-x$ if $x\in\g$; $\tau$ is called the 
\emph{principal anti-automorphism}.  

\begin{cor}\label{cor:tau_theta} 
We have $\tau(\Jla)=\Jlat$. Thus
$\tau$ induces $\theta$ according to
the commutative square:
\begin{equation}\label{eq:tau_theta} 
\begin{array}{ccc}
\Ug/\Jla&\mmaprightd{\etalas}&\Dla(X)\\[8pt]
\mapdown{\tau}&&\mapdown{\theta}\\[10pt] 
\Ug/\Jlat&\mmaprightd{\etalats}&\D^{1-\la}(X)
\end{array}  
\end{equation}
\end{cor}

\subsection{Embedding \boldmath $\Dla(X)$ into $\Dpp$}
\label{ss3_Dpp}
Let $\D(Y)$ denote  the algebra of differential operators, in the sense of
Grothendieck,  on a variety  $Y$.  
For $Y$ affine (or even quasi-affine), $\D(Y)$ coincides with the algebra
$\DD(R(Y))$ of differential operators, in the sense of
noncommutative algebra, on the  ring $R(Y)$ of regular functions.

We can identify $\pp$ with a ``big cell" $X^o$ in $X$ by means of the
Zariski open  embedding $\pp\to X$, $v\mapsto (\exp v)Q^-/Q^-$.
The canonical line bundle $\NN$ on $X$ trivializes over $X^o$; 
let $\sig$ be a   nowhere vanishing section.
Then  we have the algebra embedding  
\begin{equation}\label{eq:toDpp} 
\Dla(X)\hhookrightarrow\Dla(X^o)\mmapright{\bhla}\D(X^o)=\Dpp
\end{equation} 
where   $\bhla$ is the isomorphism defined by
$(\bhla D)(f)\sig^\la=D(f\sig^\la)$ for   $f\in R(X^o)$.
Notice that $\bhla$ is independent of the choice of $\sig$ since
$\sig$ is unique up to scaling.

We will regard (\ref{eq:toDpp}) an  inclusion. In particular
the twisted vector fields $\etala^x$ now give a realization 
of $\g$ inside $\Dpp$ where 
$\etala^x(f)=\cL_{\eta^x}(f\sig^{\la})/\sig^{\la}$.
Using the familiar rules for the Lie derivative we find
\begin{equation}\label{eq:etalax_Dpp} 
\etala^x=\eta^x+\la\left(\frac{\cL_{\eta^x}(\sig)}{\sig}\right)
\end{equation}
Thus the ``twisting" of $\eta^x$ amounts to adding a ``quantum
correction term" $\la\cL_{\eta^x}(\sig)/\sig$ which is just a 
function. We will see later that this sort of   correction is 
necessary for in order to quantize $\M$.

\begin{example}\label{ex:surr*} 
The infinitesimal action $\g\to\Vect(\pp)$, $x\mapsto\eta^x$,
integrates to a rational action of $G$ on $\pp$.
Given $Z\in\pp$, this rational action is then well-defined at $Z$ for
some  neighborhood of the identity in $G$.
In the example of \S\ref{ss2_surr}, this rational action
of $G=SL(2r,\C)$ is given by 
$\left(\Atop{A}{C}\Atop{B}{D}\right)\cdot Z=(AZ+B)(CZ+D)\inv$.
Using this, it is easy to explicitly write down the vector fields
$\eta^x$.
\end{example}

\subsection{Working in coordinates}
\label{ss3_wic} 
In this subsection we set up coordinate systems on $\pp$ and 
$\pn$. Using these we will explicitly write out realizations of
$\Ug$ in the Weyl algebras  $\Dpp$  and  $\Dpn$; see
\S\ref{ss3_etalax} and\S\ref{ss4_pilax}.

We have a unique $K$-invariant   bilinear  pairing
$\pair{\cdott}:\pp\times\pn\to\C $ such that 
$\pair{e,\ove}=r$. We can 
extend  $\pair{\cdott}$ canonically to a non-singular pairing  of the
symmetric algebras $\Spp$ and $\Spn$. Then $f\in S^d(\ppn)$ 
defines a  homogeneous degree $d$ polynomial function on $\pnp$.

Let   $v_1,\dots,v_n$ and $z_1,\dots,z_n$ be  dual vector space
bases of $\pn$  and $\pp$.  These bases form coordinate systems
on $\pp$ and  $\pn$ respectively and we can identify
\begin{equation}\nonumber 
R(\pp)=\Spn=\C[v_1,\dots,v_n]\aND
R(\pn)=\Spp=\C[z_1,\dots,z_n]
\end{equation} 
We get algebra embeddings
$S(\ppn)\to\D(\pnp)$, $A\mapsto\del_A$,
defined by $\del_{v_j}=\pdb{v_j}$ and $\del_{z_j}=\pdb{z_j}$. 
Then we can identify
\begin{equation}\nonumber 
\D(\pp)=\C[\del_{v_1},\dots,\del_{v_n},v_1,\dots,v_n]\aND
\D(\pn)=\C[\del_{z_1},\dots,\del_{z_n},z_1,\dots,z_n]
\end{equation}
We also have the  intrinsic algebra embeddings      
$S(\ppn)\to\D(\ppn)$, $A\mapsto\del^A$,
where $x\in\ppn$ defines  the constant coefficient vector field 
$\del^x$  on $\ppn$.

\subsection{The twisted vector fields  \boldmath $\eta^x_\la$}
\label{ss3_etalax} 

Let  $\nu:\q^\pm\to\C$ be the weight obtained by 
differentiating the character $\chi:Q^\pm\to\C^*$.
Recall   $r$ and  $n=\dim\ppn$ from \S\ref{ss2_tube}.
We now introduce the scalar
\begin{equation}\label{eq:nr} 
\nr=\frac{n}{r}
\end{equation}

\begin{lem}\label{lem:etalax} 
The twisted vector fields $\etala^x\in\D_1(\pp)$ are given in 
coordinates by:
\begin{equation}\nonumber
\begin{array}{lllll}
\eta^{x}_\la=-\del^x&&\mbox{\rm  if}&x\in\pp\\[7pt]
\eta^{x}_\la=-\left(\sum_{i} v_i\del^{[x,z_i]}\right)
-2\,m\,\la\,\nu(x)&&
\mbox{\rm if}&x\in\kk\\[7pt]
\eta^{x}_\la=-\half\left(\sum_{i,j} v_iv_j\del^{[[x,z_i],z_j]}\right)
+2\,m\,\la\,x
&&\mbox{\rm if}&x\in\pn
\end{array}
\end{equation}
\end{lem}
\begin{proof} 
Using the geometry of the big cell we get the coordinate expressions 
for the vector fields $\eta^x$ 
and then we work out the twisting correction
(\ref{eq:etalax_Dpp}) by choosing 
$\sig=\dd v_1\wedge\cdots\wedge\dd v_n$.
We use the fact that 
$\sig$ is $K$-semi-invariant of weight $\chi^{2\nr}$. 
See \cite{Tor},\cite{Tan}; there are minor variations in the final
answers owing to  different normalizations of $\pair{\cdott}$.
\end{proof}
We worked out these particular formulas for $\etala^x$ with Aravind
Asok in  our project on quantizing $K$-orbits in $\pn$.

So after twisting, $\pp$ acts by constant coefficient vector fields,
$\kk$ acts by homogeneous linear vector fields corrected by adding 
a constant, and $\pn$ acts by homogeneous quadratic vector fields 
corrected by adding a  homogeneous linear function.

In particular  $\pila^P$ is the constant coefficient differential   
operator   $\del^P$ if $P$ lies in $\U(\pp)=\Spp$   
(equality since   $\pp$ is abelian).

\section{Building a  noncommutative model  of   $\RM$}  
\label{sec4}

\subsection{The need for a square root}
\label{ss4_need} 

We want to try to extend our    noncommutative models
$\Dla(X)$ of $R(\OO)$ to  noncommutative models  of $\RM$.
To begin with, we observe

\begin{lem}\label{lem:C}  
Suppose $(\cC,\ga,\psi)$ is a    noncommutative  model of $\RM$. 
Extend $\psi$ to an algebra  homomorphism  $\psi:\Ug\to\cC$, 
$P\mapsto\psi^P$.  Then $\psi^F$ is a square in $\cC$ 
so that  $\psi^F=\vrho^2$ where $\vrho\in\cC_{\half[r]}$.  
\end{lem}
\begin{proof}  
By Lemma \ref{lem:BisoA},
$\cC\simeq\RM$ as   $\g$-representations with
$\cC_j\simeq\oplus_{k=0}^j\RM[k]$.
Then $\cC$ is multiplicity-free
by Proposition  \ref{prop:tO}. It follows that 
$\cC_{\half[r]}$ contains a unique copy of $V$; let $\vrho$ be a    
highest weight vector  in that copy so that 
$[\psi^x,\vrho]=\nu(x)\vrho$
for all $x\in\q^+$. Then $\vrho^2$ and $\psi^F$ are highest weight
vectors in $\cC$ of the same weight, and so they are equal up to
scaling.
\end{proof}
 
This says that we need to embed $\Dla(X)$ into some bigger algebra
where $\eta_\la^F$ becomes a square. But $\eta_\la^F=\del^F$
is a constant coefficient differential operator  and it is
uncomfortable to try to take its  square root.
It is not clear what $\sqrt{\del^F}$ could operate on. 
To remedy this,  we perform   a
Fourier transform in \S\ref{ss4_FT} following  Goncharov in \cite{Gon}. 

\subsection{Fourier transform of \boldmath $\Dla(X)$}
\label{ss4_FT}  
The Fourier transform is the   anti-isomorphism 
\begin{equation}\label{eq:FT} 
\FT:\Dpp\lto\Dpn
\end{equation}
of algebras defined by    $\FT(v)=\del^v$ for  $v\in\pn$ and 
$\FT(\del^z)=z$ for  $z\in\pp$ (see \S\ref{ss3_wic} for notations).
In \S\ref{ss3_Dpp} we embedded  $\Dla(X)$ into $\D(\pp)$. 
\begin{definition}\label{def:Bla} 
Let $\Bla=\FT(\Dla(X))$ with $\Bla_p=\FT(\Dla_p(X))$
for $p\in\bbN$. Put
$\pila^x=-\FT(\etala^{x})$ for $x\in\g$.
\end{definition}
The operators $\pila^x$ define a Lie algebra homomorphism
$\g\to\Bla_1$ which then extends to a filtered  algebra 
homomorphism
\begin{equation}\label{eq:pila} 
\pila:\Ug\lto\Bla,\qquad u\mapsto\pila^u
\end{equation}
 
\begin{lem}\label{lem:Bla_p} 
$\Bla$  is the subalgebra of $\Dpn$ generated by the   operators 
$\pila^x$ for $x\in\g$. Then $\Bla=\cup_{p\in\bbN}\,\Bla_p$ is an
algebra filtration over  $\bbN$.
The kernel of $\pila$ is $\Jlat$ and we get an induced
filtered algebra isomorphism
$\pilas:\Ug/\Jlat\lto\Bla$.
\end{lem}
\begin{proof} 
Immediate from Proposition \ref{prop:DlaX_gen} and Corollary
\ref{cor:tau_theta} since $\FT(\etala^u)=\pila^{\tau(u)}$ 
for $u\in\Ug$.  
\end{proof}

Let    $\al:\RO\to\RO$ be the Poisson algebra anti-isomorphism   
defined by $\al(\phi)=(-1)^p\phi$ if $\phi\in\RO[p]$.
Now we define $\gala$ by the   commutative square 
\begin{equation}\label{eq:ga_sq} 
\begin{array}{ccccc}
\gr\,\Dla(X)&\mmapright{\gr\FT}&\gr\Bla\\[8pt]
\mapdown{\bsla}&&\mapdown{\gala}\\[10pt] 
\RTX&\mmapright{\al\,\bj^*}&\RO
\end{array} 
\end{equation}
Proposition \ref{prop:DlaX=ncm} gives  
\begin{cor}\label{cor:Bla=ncm} 
$(\Bla,\gala,\pila)$  is a 
noncommutative model of $\RO$.
\end{cor} 

\subsection{The    operators \boldmath $\pila^x$ }
\label{ss4_pilax}  
To see what is going on, and for future use, we need to explicitly  write
out the $\pila^x$.
\begin{prop}\label{prop:pila} 
The operators $\pila^x\in\D(\pn)$ are given in coordinates by:
\begin{eqnarray}
\lefteqn{\pila^{x}=x}\hspace{9cm}& &\mbox{\rm  if }x\in\pp
\label{eq:pilax} \\[7pt]
\lefteqn{\pila^{x}=\textstyle
\big(\sum_{i} [x,z_i]\del_{z_i}\big)
+2\,m\,\la\,\nu(x)}\hspace{9cm}&&\mbox{\rm  if }x\in\kk
\label{eq:pilaxx} \\[7pt]
\lefteqn{\pila^{x}=\textstyle
\half\left(\sum_{i,j}
[[x,z_i],z_j]\,\del_{z_i}\del_{z_j}\right)
-2\,m\,\la \del^x}\hspace{9cm}&
&\mbox{\rm  if }x\in\pn\label{eq:pilaxxx} 
\end{eqnarray}
\end{prop}
\begin{proof} 
Immediate from Lemma \ref{lem:etalax}  since
$\FT(v_i)=\del_{z_i}$ and $\FT(\del_{v_i})=z_i$.
\end{proof}

So now $\pp$ acts on $R(\pn)=\Spp$ by
multiplication operators, $\kk$ acts by order $1$ differential
operators and $\pn$ acts by order $2$ differential operators.
In particular  $\pila^P=P$  if $P$ lies in $\U(\pp)=\Spp$.

\begin{cor}\label{cor:verma} 
The representation $\g\to\End\Spp$, $x\mapsto\pila^x$, is a
familiar geometric model of the generalized  Verma module for
$\q^-$ with lowest weight $\zeta=2m\la\nu$; see e.g., \cite{Tan}.
The $\g$-isomorphism is
$\Spp\to\Ug\otimes_{\U(\q^-)}\C_{\zeta}$, $u\mapsto u\otimes 1$.
The annihilator of $\Spp$ in  $\Ug$ is   $\Jlat$.
\end{cor}
 
Let $\fh\subseteq\kk$ be a Cartan subalgebra of $\kk$ and so of $\g$.
Let $\rho\in\fh^*$ be the half-sum of the positive roots
with respect to a Borel subalgebra $\fb$ of $\g$ such that 
$\fh\oplus\pp\subseteq\fb\subseteq\q^+$.
\begin{cor}\label{cor:infch} 
$\Jlat$ has infinitesimal character $-2m\la\nu+\rho$.
\end{cor}

\subsection{Extracting a square root of \boldmath $\pila^F$}
\label{ss4_sqrt} 

Our aim now is to try to extend $\Bla$ to a  noncommutative model
of $\RM$.  According to Lemma \ref{lem:C}, we need to 
extract a square root of $\pila^F$.  According to \S\ref{ss4_pilax}, this
operator  is simply the function 
\begin{equation}\label{eq:=pilaF} 
\pila^F=F 
\end{equation}

Thus we are now in a nice geometric situation, as we need to
extract a square root of the function $F$. To do this, we replace $\pn$
by its   Zariski open dense set
\begin{equation}\label{eq:Z=def}  
Z=\{q\in\pn\,|\, F(q)\neq 0\} 
\end{equation}
Then $Z$ is affine and we may identify $R(Z)=S(\pp)[F\inv]$. 
Clearly $Z$ is $K$-stable; in fact, $Z=K\cdot\ove$.
Now $F$   is not a square in
$R(Z)$; this follows for instance from Lemma \ref{lem:F}.
In the next result we construct the covering of $Z$ defined by
``extracting a square root of $F$".
\begin{lem}\label{lem:tZ} 
The complex algebraic manifold
\begin{equation}
\tZ=\{(q,t)\,|\, F(q)=t^2\}\subset Z\times\C^* 
\end{equation} 
is a non-trivial $K$-equivariant $2$-fold covering of $Z$ where the 
covering map   is $(q,t)\mapsto q$ and  $K$ acts on $\tZ$ by 
$a\cdot(q,t)=(a\cdot q,\chi(a)t)$. The formula  $w(q,t)=t$  defines a
function $w\in R(\tZ)$ such that 
\begin{equation}\label{eq:w2=F} 
w^2=F
\end{equation}
Up to isomorphism, $\tZ$ is the unique double cover of $Z$ such that
$F$ becomes a square.
\end{lem}

Notice that $\tZ$, being closed in $\Z\times\C^*$,
is affine and 
\begin{equation}\nonumber
\RtZ =\Spp{}[w\inv]=\C[z_1,\dots,z_n][w\inv]
\end{equation}
The square of the   $\C^*$-action on $Z$ lifts to the
$\C^*$-action on $\tZ$ given by $s\cdot(q,t)=(s^2q,s^rt)$.
Here  $-1$ interchanges points in  the fibers of
the cover  $\tZ\to Z$ if $r$ is odd, or acts trivially if $r$ is even.
The $\C^*$-action gives the  algebra grading 
\begin{equation}\label{eq:RtZ=oplus} 
\RtZ=\oplus_{j\in\halfZ}\,\RtZ[j]
\end{equation}
where $\RtZ[j]$ is the subspace of homogeneous functions 
of degree $2j$. Then $w$ lies in $\RtZ[{\half[r]}]$.

We have the $\g$-representation 
\begin{equation}\label{eq:Pila} 
\Pi_\la:\g\to\End\DtZ,\qquad\Pila^x(D)=[\pila^x,D] 
\end{equation}
This extends to an algebra homomorphism 
$\Pila:\Ug\to\End\DtZ$, $u\mapsto\Pila^u$.

Let $\fS$ be the Galois group of the cover $\tZ\to Z$.
Then   $\fS$ induces the algebra $\Z_2$-gradings
\begin{equation}\label{eq:grads} 
\RtZ=\RtZ\ev\oplus\RtZ\od\AND
\DtZ=\DtZ\ev\oplus\DtZ\od
\end{equation}
where $\RtZ\ev=\RtZ^{\fS}=\RZ$ and  $\DtZ\ev=\DtZ^{\fS}=\DZ$.
These gradings are  $\g$-stable in the representations
$\pila$ and  $\Pila$.

We now have    algebra inclusions 
$\Bla\subset\Dpn\subset\DZ\subset\DtZ$.
Our plan   for constructing a noncommutative model of $\RM$
is to look inside $\DtZ$  for a suitable  overring of $\Bla$.   
 
\subsection{The example \boldmath $\g=\fsl(2,\C)$}
\label{ss4_sl2} 
The simplest case occurs when $\g=\fsl(2,\C)$.
This the case $r=1$ in  \S\ref{ss2_surr} and Example \ref{ex:surr*}
and so $X=\mathbb{CP}^1$ and 
$\OO=\left\{\left(\Atop{a}{b}\Atop{c}{-a}\right)\,|\, a^2+bc=0\right\}
-\left\{\left(\Atop{0}{0}\Atop{0}{0}\right)\right\}$.
Then $\M$ identifies with $\C^2-\{0\}$  and the covering is given by 
\begin{equation}\label{eq:kap_sl2} 
\kap(\ze,\xi)=\left(\Atop{\ze\xi\phantom{o}}{\ze^2\phantom{o}}
\Atop{-\xi^2}{-\ze\xi}\right)
\end{equation}
Then $\RM=\C[\ze,\xi]$ with Poisson bracket
$\{\phi,\psi\}=\pd{\phi}{\xi}\pd{\psi}{\ze}-\pd{\phi}{\ze}\pd{\psi}{\xi}$ 
and $\RO=\C[\ze^2,\ze\xi,\xi^2]$.
The functions $\ze$ and $\xi$ are each homogeneous of degree 
$\half$  and  $\RM[j]$ is the space of homogeneous polynomials 
in $\ze$ and $\xi$ of  ordinary degree $2j$.   

Now  $v=\ove$ and $z=e$ are  dual bases of $\pn$ and $\pp$.
Twisting transforms $\eta=v^k\del_v$ into
$\etala=v^k\del_v+k\la v^{k-1}$ and so we find
\begin{equation}\nonumber 
\eta^e_\la=-\del_v,\qquad
\eta^h_\la=-2v\del_v-2\la,\qquad
\eta^{\ove}_\la=v^2\del_v+2\la v
\end{equation}   
The Fourier  transform converts these into
\begin{equation}\nonumber 
\pila^e=z,\qquad
\pila^h=2z\del_z+2\la,\qquad
\pila^{\ove}=-z\del_z^2-2\la\del_z
\end{equation}
 
We have $F=z$  and so $Z=\C^*v$ and  $w=\sqrt{z}$. 
Then $\DtZ=\C[w,w\inv,\del_w]$ where $\del_w=\pdb{w}$. 
Extending   $\DZ$ to $\DtZ$ amounts to making
the   change of variables from $z$ to $w$.
We find $\del_z=\frac{1}{2w}\del_w$ and
\begin{equation}\nonumber 
\pila^e=w^2,\qquad
\pila^h=w\del_w+2\la,\qquad
\pila^{\ove}=-\frac{1}{4}\del_w^2-
\left(\la-\tfrac{1}{4}\right)\frac{1}{w}\del_w
\end{equation}
Looking at these formulas, we see that the value
$\la=\frac{1}{4}$ is special as it eliminates the unpleasant
non-polynomial term in $\pi^e_\la$. 

So let us choose $\lao=\four$.  Then
\begin{equation}\nonumber 
\pilao^e=w^2,\qquad
\pilao^h=w\del_w+\thalf,\qquad
\pilao^{\ove}=-\four\del_w^2
\end{equation}
We recognize these operators from Weyl quantization.
They generate the even part  $\Blao=\C[w^2,w\del_w,\del_w^2]$
of the  Weyl algebra $\E=\C[w,\del_w]$. 
So $\E$ is the obvious candidate inside $\DtZ$ for an overring of $\Blao$ 
which is a noncommutative model of $\RM$. It is easy to see that
this candidate works. Indeed we  introduce the filtration
$\E=\cup_{j\in\halfN}\E_j$ where $\E_j$ is the span of the operators
$w^a\del_w^b$ for $a+b\le 2j$. 
Then $\gr\E$ is a graded Poisson algebra and we obtain the 
commutative   square
\begin{equation}\nonumber
\begin{array}{ccccc}
\gr\Blao&\hhookrightarrow&\gr\E
\\[8pt]
\mapdown{\galao}&&\mapdown{\tgalao} 
\\[10pt] 
\C[\ze^2,\ze\xi,\xi^2]&\hhookrightarrow&\C[\ze,\xi]
\end{array} 
\end{equation}
where $\tgalao$ maps the image of $w^a\del_w^b$ in 
$\gr_{\half(a+b)}\tB$ to $\ze^a\xi^b$. Now
$(\E,\tga,\pi)$ is a noncommutative model  of $\RM$ which
extends   $(\Blao,\gala,\pilao)$.
See also Remark \ref{rem:3/4}.

\subsection{The  algebras \boldmath $\Ela$}
\label{ss4_Ela}  
Our plan  is to   look inside $\DtZ$ for an extension of $\Bla$ to a 
noncommutative model $\Cla$ of $\RM$. 
Fortunately, there is a very simple way to narrow our search.
Recall that a  vector $v$ in a  $\g$-representation $\cV$ is called
\emph{$\g$-finite} if  the $\Ug$-submodule generated by $v$ is finite
dimensional. The set  $\cV\gfin$ of   $\g$-finite vectors  in $\cV$ is a 
$\g$-stable subspace which we call the \emph{$\g$-finite part} 
of $\cV$. Now   Lemma \ref{lem:BisoA} says in particular that
$\Cla$, if it exists, must lie in the $\g$-finite part of $\DtZ$.
So we make  
\begin{definition}\label{def:Ela} 
Let $\Ela$ be the $\g$-finite part of  $\DtZ$ 
in  the  representation  (\ref{eq:Pila}). 
\end{definition}
Then $\Ela$ is a subalgebra of $\DtZ$ and the action of $\fS$ defines
an algebra $\Z_2$-grading  
\begin{equation}\label{eq:Ela=ev+od} 
\Ela=\Elaev\oplus\Elaod
\end{equation}
The purpose of our next two results, Proposition \ref{prop:cap} and
Theorem \ref{thm:Ela}, is to determine the size of $\Ela$.
We will show that $\Ela$ is ``smaller than or equal to $\RM$"
in size, and moreover, if $\Cla$ exists, then $\Cla=\Ela$.

\begin{prop}\label{prop:cap} 
For each $\la\in\C$, the algebra   $\Elaev=\DZ\gfin$ is equal to $\Bla$.
\end{prop}

\begin{proof} 
Let $\FD{\ppn}$ be the fraction field of the Weyl algebra $\D(\ppn)$ 
and consider the $\g$-representations on $\FDpp$ and $\FDpn$ given
respectively by $x\mapsto[\etala^x,\cdot]$ and
$x\mapsto[\pila^x,\cdot]$.  
The Fourier transform (\ref{eq:FT}) extends uniquely to an algebra 
isomorphism $\FT:\FD{\pp}\to\FD{\pn}$.
Moreover $\FT$ identifies the $\g$-finite parts of $\FD{\pp}$ and
$\FD{\pn}$.
\begin{lem}\label{lem:gfin} 
$\Dla(X)$ is the $\g$-finite part of   $\FDpp$.
\end{lem}
\begin{proof} 
In \S\ref{ss3_Dpp} we have identified $\pp$ with a big cell in $X$ and 
the map (\ref{eq:toDpp}) embeds  $\Dla(X)$ into $\Dpp\gfin$.
To see that $\Dla(X)$ is \emph{all} of $\FDpp$, we consider
the algebra filtration $\FDpp=\cup_{j\in\Z}\,\FDpp[j]$
where $\FDpp[j]$ is the subspace spanned by the quotients
$D_1/D_2$ where $D_i\in\Dpp$ and $ord(D_1)-ord(D_2)\le j$.
(Notice that $\FDpp$ is $\Z$-filtered while $\Dla(X)$ is $\bbN$-filtered.)
Then $\gr\FDpp$ is commutative and embeds, as 
a graded Poisson algebra, into the function field $\C(\Tpp)$. 
This embedding is $\g$-linear with respect to the 
$\g$-representation on $\C(\Tpp)$ given by 
$x\mapsto\{\mu^x,\cdot\}$. 

Clearly $\RTX$ lies in $\C(\Tpp)\gfin$
and in fact $\RTX$ is \emph{all} of  $\C(\Tpp)\gfin$. To see this we
recall that,  by Proposition  \ref{prop:OO}, 
$\cN=\bj(\OO)$ is a Zariski open $G$-orbit in $T^*X$ and 
$R(\cN)=\RTX$. So 
$\C(\Tpp)\gfin=\C(T^*X)\gfin=\C(\cN)\gfin=R(\cN)=\RTX$ where
the third equality is automatic since $\cN$ is a $G$-orbit.

Thus our embedding $\Dla(X)\to\FDpp\gfin$
induces an isomorphism on the associated graded rings. In particular
then, $\gr_j\FDpp\gfin$ vanishes for $j<0$ and so 
$\FDpp[j]\gfin=\FDpp[-1]\gfin$ if $j<0$. But $\cap_{j<0}\,\FDpp[j]=0$  
and so  $\FDpp[j]\gfin=0$ if $j<0$. It follows now that
$\Dla_j(X)=\FDpp[j]\gfin$ for all $j\in\Z$.
\end{proof}

Thus $\FDpn\gfin=\FT(\FDpp\gfin)=\FT(\Dla(X))=\Bla$.
But also $\FDpn\supset\DZ\supset\Bla$. So
$\DZ\gfin=\Bla$.
\end{proof}

\section{The algebras $\Ela$ and symplectic geometry of $\M$}  
\label{secE}

\subsection{Filtration theorem for \boldmath $\Ela$} 
\label{ssE_thm}

We  regard our $\bbN$-filtration of $\Bla$ as a 
$\halfN$-filtration by the recipe given in \S\ref{ss3_ncm}.
Recall $\Ela=\Bla\oplus\Elaod$ by Proposition \ref{prop:cap}.

\begin{thm}\label{thm:Ela} 
Pick $\la\in\C$. 
There is a unique $\fS$-stable algebra filtration 
\begin{equation}\label{eq:Ela=cup} 
\Ela=\cup_{p\in\halfN}\,\Ela_p
\end{equation} 
extending our   filtration of $\Bla$ such that $\gr\Ela$ has no 
zero-divisors. This satisfies
$[\Ela_p,\Ela_q]\subseteq\Ela_{p+q-1}$ and so $\gr\Ela$ is a
graded Poisson algebra.

The map $\gala:\gr\Bla\to\RO$  defined in \textup{(\ref{eq:ga_sq})}
extends, uniquely up to  the $\Z_2$-actions defined by 
the Galois groups $\fS$ and $\cS$,
to a $\Z_2$-equivariant Poisson algebra homomorphism  
\begin{equation}\label{eq:tga} 
\tgala:\gr\Ela\lto\RM 
\end{equation}
In fact  $\tgala$ is $1$-to-$1$. Thus, if we identify $\gr\Ela$ with its
image, we get
\begin{equation}\label{eq:subsub} 
\RO\subseteq\gr\Ela\subseteq\RM
\end{equation}
\end{thm}
\S\ref{ssE_proof} is devoted to proving  Theorem \ref{thm:Ela}.
  
\subsection{Proof of Theorem \ref{thm:Ela}}
\label{ssE_proof}
Suppose we have extended our filtration of $\Bla$ to an algebra
filtration (\ref{eq:Ela=cup}) such that $\gr\Ela$ has no nilpotents.
Say $S\in\Elaod$. Then $S^2$ lies in $\Bla$ and so has 
some known filtration  degree $p$. 
Since  $\gr\Ela$ has no zero divisors, it follows that
the filtration degree of a product is equal to the sum of the 
filtration degrees of the factors. Hence
$S$ has filtration degree $\half p$. This proves uniqueness.

To prove existence, we will construct an algebra  filtration of
$\DtZ$ and then restrict it to $\Ela$.   We start with 
the  vector space  isomorphism
\begin{equation}\label{eq:bm} 
\bm:\RTtZ=\RtZ\otimes\Spn\lto\DtZ,\qquad\bm(f\otimes P)=f\del^P 
\end{equation}
Here  $\del^P\in\Dpn$ defines  a differential operator on $\tZ$
by first restricting it to $Z$ and then lifting it to $\tZ$. Then
$P$ defines a function on $\TtZ$, namely the symbol of $\del^P$. 

For $j\in\halfZ$,  we put  
\begin{equation}\nonumber 
\DtZs[j]=\bm\left(\RtZ[\le j]\otimes\Spn\right) 
\end{equation}
Then $\DtZ=\cup_{j\in\halfZ}\,\DtZs[j]$ is an $\fS$-stable 
algebra  filtration. This induces a filtration
$\DZ=\cup_{j\in\halfZ}\,\DZs[j]$  which extends  our  filtration 
of $\Bla$.  We have $[\DtZs[j],\DtZs[k]]\subseteq\DtZs[j+k-1]$.  
Put $\grs_j\DtZ=\DtZs[j]/\DtZs[j-\half]$. 
Let $\bps_j:\DtZs[j]\to\grs_j\DtZ$ be the natural projection. Then 
$\grs\DtZ=\oplus_{j\in\halfZ}\grs_j\DtZ$
is commutative and acquires a Poisson 
bracket induced by commutator of operators.  

In this way,
$\grs\DtZ$ is a graded Poisson algebra in the sense of 
Definition \ref{def:gPa}  with $\halfZ$ replacing $\halfN$.
Although this $\sharp$-filtration of $\DtZ$ is not the usual (order)  
filtration,  fortunately $\grs\DtZ$ is very nice. It is easy to check  
\begin{lem}\label{lem:btt} 
The   map $\bm\inv:\DtZ\to\RTtZ$ 
induces an $\fS$-invariant  Poisson algebra isomorphism
$\btt:\grs\DtZ\to\RTtZ$. 
\end{lem}

For $D\in\DtZs[j]$, we say that $\btt(\bps_j(D))\in\RTtZ$ is the   
\emph{degree $j$ $\sharp$-symbol of $D$}.
What is happening here is that the two associated graded algebras
of $\DtZ$, for the $\sharp$-filtration and  the order filtration,
are the same as Poisson algebras but of course different as 
graded algebras.  We call the  grading  
$\RTtZ=\oplus_{j\in\halfZ}\,\btt(\grs_j\DtZ)$ the 
\emph{$\sharp$-gradation}; this corresponds  to homogeneity
along the base $\tZ$.

\begin{example}\label{ex:sl2_sharp} 
In Example \ref{ss4_sl2},  
$w^a\del_z^b=z^{\half[a]}z^{\half[b]}\del_z^b$ has
$\sharp$-filtration degree  equal to $\half(a+b)$. 
\end{example}

For  $j\in\halfZ$, we put
\begin{equation}\label{eq:Elaj=} 
\Ela_j=\Ela\cap\DtZs[j]
\end{equation}
Then $\gr\Ela$ is a Poisson subalgebra of $\grs\DtZ$ and this  has no zero 
divisors by  Lemma \ref{lem:btt}. 
We see now that $\Ela=\cup_{j\in\Z}\,\Ela_j$ is an algebra filtration
which has all the desired  properties,  
except  that we still need to prove  
$\Ela_j=0$ for $j<0$. This   will be easy once we  analyze $\gr\Ela$.

The $\g$-representation (\ref{eq:Pila}) induces a $\g$-representation 
on $\RTtZ$ and this is given by  $x\mapsto\{b^x,\cdot\}$ where  
$b^x\in\RTZ$ is the degree $1$ $\sharp$-symbol of $\pila^x\in\DZs[1]$. 
Now  $\gr\Ela$  lies inside $\RTtZ\gfin$ 
(recall Definition \ref{def:Ela}) and we get the commutative diagram
\begin{equation}\label{eq:grEla_rect} 
\begin{array}{ccccccc}
\gr\Bla&\hhookrightarrow&\gr\Ela&\hhookrightarrow&\grs\DtZ
\\[8pt]
\mapdown{\btt}&&\mapdown{\btt}&&\mapdown{\btt} 
\\[10pt] 
\RTZ\gfin&\hhookrightarrow&\RTtZ\gfin&\hhookrightarrow&\RTtZ 
\end{array} 
\end{equation}

The problem now is to recognize
the functions $b^x$ and then compute  $\RTtZ\gfin$. Let us identify
$R(T^*\ppn)=\Smp\otimes\Spm$ and let
$\bx: \RTpp\lto\RTpn$ be the linear map such that 
$\bx(f\otimes g)=g\otimes f$.

\begin{lem}\label{lem:nux} 
For each $x\in\g$, $b^x$ lies in $\RTpn$ and 
is the image under $\bx$ of the usual order $1$
symbol of the twisted vector field $\etala^{-x}$ on $\pp$.
\end{lem}
\begin{proof} 
Using Lemma \ref{lem:btt} we find the commutative diagram
\begin{equation}\label{eq:grsDZ_rect} 
\begin{array}{ccccccc}
\gr\,\Dla(X)&\hookrightarrow&
\gr\Dpp&\mmaprightd{\gr\FT}&\grs\Dpn&\hookrightarrow&\grs\DZ
\\[8pt]
\mapdown{\bsla}&&\mapdown{\bsla}&&\mapdown{\btt}&
&\mapdown{\btt} 
\\[10pt] 
\RTX&\hookrightarrow&
\RTpp&\mmaprightd{\bx}&\RTpn&\hookrightarrow&\RTZ
\end{array} 
\end{equation}
Then commutativity of the middle square gives 
$\bx(\bsla(\etala^x))=\btt(\pila^{-x})=b^{-x}$.
\end{proof}

Let $\bfb:\TZ\to\g$ be the moment map defined by 
$\Killg{\bfb(m),x}=b^x(m)$.
\begin{prop}\label{prop:tbfb} 
The moment map $\bfb$ is a  symplectic Zariski open embedding
of $\TZ$ into $\OO$. Moreover $\bfb$ lifts to
a   $\Z_2$-equivariant symplectic Zariski open embedding $\tbfb$
of $\TtZ$ into $\M$ so that we get the commutative square
\begin{equation}\label{eq:TtZ_sq} 
\begin{array}{ccccc}
\TtZ&\mmapright{\tbfb}&\M
\\[8pt]
\mapdown{}&&\mapdown{\kap} 
\\[10pt] 
\TZ&\mmapright{\bfb}&\OO
\end{array} 
\end{equation}
\end{prop}
\begin{proof} 
To prove this is it easiest to  start with the geometry of $\OO$ and $\M$.
In \S\ref{ss4_sqrt}, we introduced $Z$ and $\tZ$ expressly for the
purpose of extracting a square root of the function $F$.
In fact, $Z$ occurs naturally in the geometry of $\OO$. To begin with,  
$Z=K\cdot\ove=\OO\cap\pn$ and so $Z$ is a smooth Lagrangian 
submanifold of $\OO$ with respect to the  KKS symplectic form $\om$.
Then   $F\in\RZ$ is simply  the restriction of $\phiF\in\RO$, i.e.,
\begin{equation}\label{eq:phiF=F} 
\phiF|_Z=F
\end{equation}

The composite map $\OO\mapright{\bj}\,T^*X\to X$ makes $\OO$ into a  
$G$-equivariant fiber bundle  over $X$  with typical fiber
$Z$.  Indeed, the cotangent bundle $T^*X\to X$
identifies with the contracted product bundle
$G\times_{Q^-}\pn\to G/Q^-=X$ and then $\OO$ identifies with
$G\times_{Q^-} Z$.   We will  treat the map $\bj$ as an inclusion.

The cotangent bundle $T^*X\to X$, and hence the subbundle
$\OO\to X$, trivializes over the big cell $X^o$. We have identified
$X^o$ with $\pp$.  Now  we get
the following commutative diagram: 
\begin{equation}\label{eq:rect_FT} 
\begin{array}{ccccccccc}
T^*X&\supset&\Tpp&=&\pp\times\pn&\xrightarrow{permute}&
\pn\times\pp&=&\Tpn\\[5pt]
\mapup{\bj}&&\mapup{}&&\mapup{}&&\mapup{}&
&\mapup{}\\[9pt] 
\OO&\supset&\OO\cap\Tpp&=&\pp\times Z&
\xrightarrow{permute}& Z\times\pp&=&T^*Z 
\end{array} 
\end{equation}
Here all maps are birational symplectomorphisms, 
except the two permutation maps are anti-symplectic.
The mathematical content of the left part of (\ref{eq:rect_FT}) 
is that $\OO\cap\Tpp$ is a trivial bundle over $\pp$
and moreover the standard trivialization
$\Tpp=\pp\times\pn$ induces the trivialization 
$\OO\cap\Tpp=\pp\times Z$. 
This is true  because $\pp$ is abelian.

The bottom row of (\ref{eq:rect_FT})  read right to left
defines an anti-symplectic Zariski open embedding  $\TZ\to\OO$. 
Using Lemma \ref{lem:nux} we see that the composition of this   
embedding with the  map $\OO\to\g$, $w\mapsto -w$, 
is the  moment map $\bfb:\TZ\to\g$. This proves the first statement. 
 
The  covering $\kap:\M\to\OO$  induces a covering $\kap\inv(Z)\to Z$.
This identifies with the  covering
$\tZ\to Z$ we constructed in  Lemma \ref{lem:tZ} (on account of 
(\ref{eq:phiF=F}) for instance) in such a way that
we get the   commutative square 
\begin{equation}\label{eq:tZ_sq_intro} 
\begin{array}{ccccc}
\tZ&\hhookrightarrow&\M
\\[8pt]
\mapdown{}&&\mapdown{\kap} 
\\[10pt] 
Z&\hhookrightarrow&\OO
\end{array} 
\end{equation}
where $\ze|_{\tZ}=w$  (cf. (\ref{eq:ze=}) and (\ref{eq:w2=F})). Then
the inclusion of $\tZ$ into $\M$ is $K$-invariant,
the Galois groups $\cS$ and $\fS$ identify naturally and 
$\tZ$ is Lagrangian in $\M$.

Now let $\cN=\kap\inv(\TZ)$.   We can lift the projection  
$\tau:\cN\;\mapright{\kap}\;\TZ\to Z$  to a map 
$\ttau:\cN\to\tZ$  in the following way.
Notice that if $p,p'\in\tZ$ lie above $q\in Z$, then
$\kap\inv(T^*_qZ)$ breaks into two connected components
$\cN_{p}$ and $\cN_{p'}$ which contain $p$ and $p'$ respectively.
Then $\cN=\cup_{p\in\tZ}\cN_p$. Now we define
$\ttau$ by $\ttau(\cN_p)=p$.
Then $\cN$  identifies naturally with 
$\TZ\times_Z\tZ\simeq\TtZ$ and the rest of the result follows.
\end{proof}
\begin{cor}\label{cor:RTtZ} 
We have the  commutative diagram 
\begin{equation}\label{eq:RTtZ_sq} 
\begin{array}{ccccc}
\RTtZ\gfin&\mmapleft{\tbfb^*}&\RM
\\[8pt]
\mapup{}&&\mapup{\kap^*} 
\\[10pt] 
\RTZ\gfin&\mmapleft{\bfb^*}&\RO
\end{array} 
\end{equation}
The  horizontal maps are $\g$-linear
graded Poisson algebra isomorphisms, where $\RTtZ$ and $\RTZ$ 
have the  $\sharp$-gradations.
\end{cor}
\begin{proof} 
The only point that is not immediate is the surjectivity of
$\tbfb^*$ and $\bfb^*$ in (\ref{eq:RTtZ_sq}). But this follows
since, as  $\M$ and $\OO$ are   $G$-orbits, 
$\RM$ and $\RO$ are  the $\g$-finite parts of the function fields 
$\C(\M)$ and $\C(\OO)$.
\end{proof}

Corollary \ref{cor:RTtZ} says in particular  that the 
$\sharp$-gradation of $\RTtZ\gfin$, and hence our gradation
of $\gr\Ela$,   vanishes in negative degrees.
So $\Ela_{j}=\Ela_{-1/2}$ if $j<0$.  But $\cap_{j<0}\Ela_{j}=0$ since
each operator $D\in\DtZ$ has finite $\sharp$-filtration degree.
Thus $\Ela_{j}=0$ for all $j<0$.

Let $\tgala$ be the composite map 
\begin{equation}\label{tga=}  
\gr\Ela\mmapright{\btt}\RTtZ\gfin\xrightarrow{(\tbfb^*)\inv}\RM
\end{equation}
Clearly $\tgala$ extends $\gala$ and $\tgala$ is a $\Z_2$-equivariant 
embedding of Poisson algebras. 

Finally   suppose that $\al:\gr\Ela\to\RM$ is some map enjoying
the same properties as $\tga=\tgala$.  
We want to to show that $\al$ is either $\tga$ or  
$\vsig\tga$ where $\vsig$ is the non-trivial automorphism of $\RM$
defined by the $\cS$-action. 
Since $\al$ is $1$-to-$1$ on $\gr\Bla$ it follows easily  that
$\al$ is $1$-to-$1$ on $\gr\Ela$.
Then  by considering  the fraction field   of $\al(\gr\Ela)$  
we find  $\al=\tga$ or $\al=\vsig\tga$.

This completes the proof of  Theorem \ref{thm:Ela}.  

\begin{remark}\label{rem:tbfb}  
(i)  If  $x$ lies in one of $\kk,\pp,\pn$, 
then $b^x$  is the symbol (in the usual sense)  of $\pila^x$.
(ii) the map  $\tbfb$  in Proposition \ref{prop:tbfb} 
is   unique up to the $\Z_2$-action.
\end{remark}

\subsection{Simplicity of \boldmath  $\Ela$}
A nice fact which will be important later   (see Corollary
\ref{cor:Ela0=simple} and Proposition \ref{prop:T}) is 
\label{ssE_simple}
\begin{cor}\label{cor:Ela=simple} 
$\Ela$ is a simple ring if $\Dla(X)$ is a simple ring.
\end{cor}
\begin{proof} 
Let $L$ be a non-zero $2$-sided ideal in $\Ela$.
Then $L$ is in particular $\g$-stable with respect to the 
representation (\ref{eq:Pila}). Let $L_0\subset L$ be a 
subspace carrying a non-zero $\g$-irreducible representation.
Then $L_0$ lies in  $\Elaev$ or  $\Elaod$. This follows since, 
by   (\ref{eq:subsub}), 
$\Ela$ is isomorphic as a $G$-representation to a subspace of
$\RM$ and so, by Proposition \ref{prop:tO}, $\Ela$ is multiplicity-free.
So $L_0^2$ lies in $\Elaev$ which is equal  to $\Bla$  by Proposition
\ref{prop:cap}.  So $L\cap\Bla\neq 0$.
But $\Bla$ is anti-isomorphic to $\Dla(X)$ and so is simple by
hypothesis. Thus  $L$ contains $\Bla$ and so $L$ contains $1$.
\end{proof}

\subsection{The anti-symmetry \boldmath $\la\mapsto(1-\la)$}
\label{ssE_1-la}
We define
an \emph{anti-isomorphism of noncommutative models} from
$(\Ela,\tgala,\pila)$ to   $(\E^{\la'},\tga_{\la'},\pi_{\la'})$
to be  a filtered algebra anti-isomorphism $\delta:\Ela\to\E^{\la'}$
such that we have $\delta(\pila^x)=-\,\pi_{\la'}^x$ and  
commutativity in
\begin{equation}\label{eq:grdelta_sq} 
\begin{array}{ccccc}
\gr\,\Ela &\mmapright{\gr\delta}&\gr\E^{\la'}\\[8pt]
\mapdown{\tgala}&&\mapdown{\tga_{\la'}}\\[10pt] 
\RM&\mmapright{\al}&\RM
\end{array} 
\end{equation}
where $\al$ is  the automorphism   of $\RM$ 
defined by  $\al(\phi)=i^{2j}\phi$ if $\phi\in\RM[j]$. 
(So we are  extending the involution $\al$ of $\RO$ defined before
(\ref{eq:ga_sq}).)   Notice that then $\delta$ is $\g$-linear, and so
$G$-equivariant, with respect to the representations $\Pila$ and $\Pilap$.
Now  Proposition \ref{prop:theta} gives
\begin{cor}\label{cor:tdelta} 
The anti-isomorphism   $\FT\theta\FT\inv$ extends,  
uniquely up the action of  $\fS$,  
to   a  $(\GfS)$-invariant anti-isomorphism $\delta$ so that we get 
the   commutative diagram 
\begin{equation}\label{eq:sig_rect} 
\begin{array}{ccccccc}
\Dla(X)&\mmaprightd{\FT }&\Bla&\hhookrightarrow{}&\Ela\\[8pt]
\mapdown{\theta}&&\mapdown{\FT\theta\FT\inv}
    &&\mapdown{\delta}\\[8pt]
\D^{1-\la}&\mmaprightd{\FT}&\B^{1-\la}
&\hhookrightarrow{}&\E^{1-\la}
\end{array}   
\end{equation}
We can specify $\delta$ by $\delta(w)=i^rw$ and then
$\delta$ gives 
an $\fS$-invariant anti-isomorphism of noncommutative models
from $(\Ela,\tgala,\pila)$ to $(\E^{1-\la},\tga_{1-\la},\pi_{1-\la})$. 
\end{cor}
\begin{proof} 
Define  $\delta:\Bla\to\B^{1-\la}$ by $\delta=\FT\theta\FT\inv$.
Then $\delta$   extends to an algebra  anti-involution
$\delta:\Dpn\to\Dpn$ 
where $\delta(z_i)=-z_i$ and $\delta(\del_{z_i})=\del_{z_i}$.
This follows from the definition of $\theta$ in the proof of 
Proposition \ref{prop:theta}.
Now $\delta$ naturally extends  to an anti-involution of $\DZ$ and 
then  to an $\fS$-invariant  anti-automorphism $\delta$ of     
$\DtZ$ such that $\delta(w)=i^rw$.  Then $\delta$ preserves the
$\sharp$-filtration of $\DtZ$. The relation
$\theta(\etala^x)=-\,\eta_{1-\la}^x$
implies $\delta(\pila^x)=-\,\pi_{1-\la}^x$.
So  $\delta(\Ela)=\E^{1-\la}$. The rest  is now clear. 
\end{proof}

\subsection{\boldmath $\g$-finiteness of $w$ }
\label{ssE_gfin}
Let $\Wla=\Ug\cdot w$ be the $\Ug$-submodule of $\DtZ$ generated by
$w$ in the representation (\ref{eq:Pila}). So
$w$ is $\g$-finite $\Leftrightarrow$ $\Wla$ is finite-dimensional.
Theorem \ref{thm:Ela} gives 
\begin{cor}\label{cor:Ela} 
Pick $\la\in\C$. The following are equivalent:
\begin{itemize}
\item[\rm(i)] $w\in\Ela$, i.e., $w$ is $\g$-finite
in the representation \textup{(\ref{eq:Pila})}.
\item[\rm(ii)] $\tgala(\gr\Ela)=\RM$.
\item[\rm(iii)] The triple $(\Ela,\tgala,\pila)$ is a  noncommutative
model of $\RM$.
\item[\rm (iv)] $\Ela$ is generated as an algebra by $\Bla$ and $w$.
\item[\rm (v)]  We have $\Ela=\Bla\oplus\Wla\Bla$.
\end{itemize}
\end{cor}
\begin{proof} 
The equivalence (ii)$\Leftrightarrow$(iii) is immediate from 
Theorem \ref{thm:Ela}. The implication (i)$\Rightarrow$(ii) 
follows easily from (\ref{eq:RtO=}). We get
(iii)$\Rightarrow$(i) as follows. Given (iii),  we know by Lemma
\ref{lem:C} that $\pila^F=F$ admits  a square root $D\in\Ela$.
But the equation  
$D^2=F$ in $\DtZ$ forces $D\in\DtZ[0]=\RtZ$. 
So $D=\pm w$. So $w\in\Ela$.
Finally, the equivalences  with (iv) and (v)  
follow  from Corollary \ref{cor:RtO=PA} and (\ref{eq:RtO=}).
\end{proof}
The problem now is to determine which, if any, values of $\la$
satisfy (i)-(v); we call these critical values.  We solve this in Theorem
\ref{thm:la} below.

\section{Critical values of   $\la$}  
\label{sec6}
\subsection{Critical values theorem}
We find exactly two critical values of $\la$. Recall the number $m$
defined by (\ref{eq:nr}) and Table \ref{table:O}.
\label{ss6_critical} 
\begin{thm}\label{thm:la} 
There are exactly two values of $\la$, namely
$\la=\thalf\pm\frac{1}{4\nr}$,
such that the triple $(\Ela,\tgala,\pila)$ constructed in 
\textup{\S\ref{sec4}-\S\ref{secE}} is a noncommutative model of    
 $\RM$.
\end{thm}
The proof occupies  \S\ref{ss6_S}-\S\ref{ss6_com}.  
From now on we write
\begin{equation}  
\lao=\thalf-\frac{1}{4\nr}\AND \laop=\half+\frac{1}{4\nr}
\end{equation}

\begin{remark}\label{rem:thmla}  
The two algebras $\Elao$ and $\Elaop$ are anti-isomorphic by 
Corollary \ref{cor:tdelta}  since $\lao+\laop=1$.
The difference
$d=\lao-\laop=-\frac{1}{2\nr}$ has the property that 
$\NN^{-d}$ is the  $G$-homogeneous line bundle over $X$ with
$\Gamma(X,\NN^{-d})\simeq V$.  With a little thought we see   that
this is what we would  expect. 
It might interesting to interpret the fact that $\lao$ and $\laop$ tend to
$\half$ as the rank of $\g$ tends to infinity.
\end{remark}

Here is an overview of the proof of Theorem \ref{thm:la}.
By Corollary \ref{cor:Ela},   Theorem \ref{thm:la}  amounts to 
\begin{equation}\label{eq:wiff} 
\mbox{$w$ is  $\g$-finite in $\DtZ$} \quad\Leftrightarrow\quad
\la=\thalf\pm\frac{1}{4\nr}
\end{equation}
where   we are considering
$\g$-finiteness of $w$ with respect to  the representation 
(\ref{eq:Pila})  which of course  depends on $\la$. 

In  \S\ref{ss6_S}, we reduce proving  (\ref{eq:wiff}) to
verifying that  a certain double commutator  vanishes in $\DtZ$. In
\S\ref{ss6_com} we carry out the computations. To do this,
we exploit the fact that
$\pn$  is a Jordan algebra and $Z\subset\pn$ is the subset of
Jordan invertible elements.  A key point is that the Jordan
theory gives us explicit formulas for the first
and second partial derivatives of $w$ and we explain this in
\S\ref{ss6_Jordan}.
To aid the reader, we   give several
explicit references to  \cite{F-K}; see also  \cite{Sat}. 

We are using Jordan algebra theory as a tool.  It would be very
interesting to find some deeper connections, which we feel surely
exist. 

\begin{remark}\label{rem:3/4}  
In the example of \S\ref{ss4_sl2}, 
Theorem \ref{thm:la} gives the two critical values  
$\la=\four,\frac{3}{4}$. In hindsight, we see that the ``new" value
$\la=\frac{3}{4}$  arises by
writing our differential operators in the form
$\del_w^aw^b$. Indeed, then $\la=\frac{3}{4}$ becomes the value for
which an unpleasant term goes away.
\end{remark}

\subsection{Strategy }
\label{ss6_S} 
In the representation  (\ref{eq:Pila}) we have 
$[\pila^x,w]=\nu(x)w$ for all $x\in\q^+$ and so
$\Wla=\Spn\cdot w$. Similarly, in  the representation (\ref{eq:Phi}),
we have $\{\phi^x,\ze\}=\nu(x)\ze$ for all $x\in\q^+$ and 
so $V=\Spn\cdot\ze$. 
Let $\Ila\subset\Spn$ and $\I\subset\Spn$ be the annihilators of  
$w\in\Wla$ and $\ze\in V$   respectively. Since $V$ is finite-dimensional,
It follows by highest weight theory that, as a $\Ug$-module,
$\Wla$ admits a unique  finite-dimensional quotient and this is 
isomorphic to $V$.  Hence $\Ila\subseteq\I$ and $\Wla$ is
finite-dimensional  if and only if  $\Ila=\I$, i.e.,
\begin{equation}\label{eq:iff} 
\mbox{$w$ is $\g$-finite}\quad\Leftrightarrow\quad \Ila=\I
\end{equation}

Using the action of $\adj h$, we see that $\I$ and $\Ila$ are graded
ideals   in $\Spn$ so that  $\I=\oplus_{p\ge 0}\I^p$ and 
$\Ila=\oplus_{p\ge 0}\Ila^p$.
We can give  a precise description of $\I$ once we
recall the structure of $\Spn$ as a $K$-representation. 
 
To do this, we need to set up some structure to write down
lowest weights of $K$-representations. 
We can find $r$ commuting Lie subalgebras $\s_1,\dots,\s_r$ of $\g$,
such that (i) each $\s_i$ is isomorphic to $\fsl(2,\C)$, 
(ii) the sum $\s_1+\cdots+\s_r$ is direct and contains 
$\s=\C e\oplus\C h\oplus\C\ove$,
(iii) the   decompositions 
$e=\sum_{i=1}^re_i$ and $h=\sum_{i=1}^rh_i$, 
where $e_i,h_i\in\s_i$, satisfy $e_i\in\pp$ and $h_i\in\kk$, and
(iv) each $\s_i$ is stable under complex conjugation and 
under the complex Cartan involution of  $(\g,\kk)$.   
Then $h_1,\dots,h_r$ span an $r$-dimensional abelian subalgebra
$\fa\subset\kk$ and  $e_i$ is a weight vector of $\fa$  of weight
$2t_i$ where $t_i(\sum_{p=1}^rc_ph_p)=c_i$.
 
The pair $(\kk,\ks)$ is a complex symmetric pair; we have  
$\ks=\kk^e$. We can embed $\fa$ in a Cartan subalgebra
$\fh\subset\kk$ such that $\fh=\fa\oplus(\fh\cap\ks)$. 
Now we extend  $t_1,\dots,t_r$ to  weights of $\fh$ 
by having them vanish on $\fh\cap\kk^e$.
We set $\ga_i=-2(t_1+\dots+t_i)$.
The weights $\ga_1,\dots,\ga_r$ are (for
an appropriate system of positive roots) the lowest weights of $r$ 
distinct finite-dimensional irreducible $\kk$-representations. In
particular,  $-2\nu=\ga_r$.
 
\begin{thm}\label{thm:Sch}\textup{\cite{Sch}}
$\Spn$ is a multiplicity-free $K$-representation and the set of lowest
weights occurring in $\Spn$ is
$\{\sum_{i=1}^rc_i\ga_i\,|\, c_i\ge 0\}$.  We have
\begin{equation}\label{eq:Sqpn=} 
S^q(\pn)=\bigoplus_{c_1+2c_2+\cdots+rc_r=q}
L_{c_1\ga_1+\cdots+c_r\ga_r}
\end{equation}
where the subspace $L_{\psi}$ carries the  $K$-representation 
of lowest  weight   $\psi$. 
\end{thm}

\begin{lem}\label{lem:I2} 
The ideal $\I\subset\Spn$ is generated by its degree $2$ component
$\I^2$ and $\I^2=L_{2\ga_1}$.  Then: 
$w$ is $\g$-finite $\Leftrightarrow$ $L_{2\ga_1}\subset\Ila$. 
 \end{lem} 
\begin{proof} 
By decomposing $V\simeq \Spn/\I$ as a $K$-representation, we find
that the (unique) $K$-stable direct sum complement to $\I$ in $\Spn$ is
$\oplus_{i=0}^r\,L_{\ga_i}$.
So in particular, $L_{\ga_2}$ is the complement to $\I^2$ in
$S^2(\pn)$. But $S^2(\pn)=L_{2\ga_1}\oplus L_{\ga_2}$  
and so $\I^2=L_{2\ga_1}$. We find that $\I^2$ generates  $\I$, 
and so $\Ila=\I$ $\Leftrightarrow$ $\I^2\subset\Ila$.
\end{proof}

We can simplify the criterion in Lemma \ref{lem:I2} considerably.
\begin{lem}\label{lem:wgfin} 
Let    $y=\ove_1$. Then:  $w$ is $\g$-finite $\Leftrightarrow$   
$[\pila^y,[\pila^y,w]]=0$.
\end{lem}
\begin{proof} 
Since  $u=\ove_1^2$ is a lowest 
weight vector in $L_{2\ga_1}$ and $w$ is $K$-semi-invariant, it   
follows that $L_{2\ga_1}\subset\Ila$ $\Leftrightarrow$ 
$\Pila^u(w)=0$. But $\Pila^u(w)=[\pila^y,[\pila^y,w]]$.
\end{proof}

\subsection{Calculus on the coupled Jordan algebras \boldmath
$\ppn$}  
\label{ss6_Jordan} 
By TKK theory, $\pp$ and $\pn$ are (isomorphic) complex simple    
Jordan algebras with  Jordan products defined   by
$[x,e]\circ [y,e]=[x,[y,e]]$ and 
$[x,\ove]\circ [y,\ove]=[x,[y,\ove]]$ where $x,y\in\fr$
and $\fr$ is the orthogonal complement  in $\kk$ to $\ks$.
The Jordan identity elements are $e\in\pp$ and $\ove\in\pn$.
In this subsection, we explain some basic formulas from Jordan theory
that we will use throughout \S\ref{ss6_com}.
See \cite[Table on page 160]{F-K} for the list of complex Jordan algebras
carried by $\ppn$.  (Notation warning:  our algebra $\kk$ is called  ``$\g$" 
in  \cite{F-K}.)

From now on, we usually omit the symbol
"$\circ$" and write  the Jordan product $a\circ b$ as  $ab$.
We adopt this convention: if $D$ is an operator and $f$ is a function
then $Df$ is the composition of operators (where $f$ is regarded as a
multiplication operator) and $[Df]$ is the function obtained by
applying $D$ to $f$.

The polynomial function $F$, normalized so that $F(\ove)=1$,
is the Jordan norm of $\pn$.
Thus  by (\ref{eq:Z=def}), $Z$ is the set of Jordan invertible elements    
in $\pn$. The partial derivatives of $F$ are
(see \cite[Prop. III.4.2, page 52]{F-K}), where  $v\in\pn$,
\begin{equation}\label{eq:delvF} 
[\del^{v}F]=F\tr(v\,q\inv)
\end{equation}
Here  $q$ is an arbitrary  point in $\pn$ so that the RHS of 
(\ref{eq:delvF})
is the function  $q\mapsto F(q)\tr(v\,q\inv)$ where $q\inv$ is the
Jordan inverse and  $\tr$ is the Jordan trace.
A quick definition is  $\tr(x)=\frac{1}{m}\Tr L_x$ where 
$L_x:\ppn\to\ppn$ is   Jordan multiplication by
$x\in\ppn$ and  $\Tr$ is the usual trace. Then 
$\tr(e)=\tr(\ove)=r$.  See  \cite[II.2 and Prop. III.4.2]{F-K}. 
 
Since $w=\sqrt{F}$, (\ref{eq:delvF}) gives 
\begin{equation}\label{eq:delvw} 
[\del^v{w}]=\half w\tr(v\,q\inv)
\end{equation}
To get  $[\del^u\del^vw]$, where $v,u\in\pn$,  we recall
(\cite[Proposition II.3.3(i), page 33]{F-K})  
\begin{equation}\label{eq:delutr} 
[\del^u\tr(vq\inv)]=-\tr(v\jtp{q\inv,u,q\inv})
\end{equation}
Here $\jtp{a,b,c}=a(bc)+b(ac)-(ab)c$ is the 
\emph{Jordan triple  product}.    Then
\begin{equation}\label{eq:deldelw} 
[\del^{u}\del^{v}w] 
=\four w\tr(u\,q\inv)\tr(v\,q\inv)-\half w
\tr(v\jtp{q\inv,u,q\inv})
\end{equation}

The coupling of our Jordan algebras is achieved by the transpose maps
$\ppn\to\pnp$, $x\mapsto x^t$, defined by
$\pair{x,y}=\tr(xy^t)=\tr(x^ty)$. These maps are 
inverse Jordan algebra isomorphisms and $\tr(u)=\tr(u^t)$.
From now on, we assume that our basis $v_1,\dots,v_n$ of
$\pn$ introduced in \S\ref{ss3_wic}  is
orthonormal with respect to $\tr$. Then  $v_i^t=z_i$
and $z=\tr(z^tq)$.

Our realization $x\mapsto\eta^x$ of $\g$ inside
$\Dpp$ is the  TKK construction (see \cite{Sat}). In  that
language   (\ref{eq:pilaxxx}) becomes, for  $y\in\pn$,
\begin{eqnarray}
\pila^{y}&=&\textstyle -\left(\sum_{i,j}\jtp{z_i,y^t,z_j}
\,\del_{z_i}\del_{z_j}\right)-2\,m\la \del^y\nonumber\\[3pt]
&=&\textstyle -\left(\sum_{i,j}\tr(\jtp{v_j,y,v_i}q)
\,\del^{v_i}\del^{v_j}\right)-2\,m\la \del^y\label{eq:pilaxxx*} 
\end{eqnarray} 
since $\jtp{a,b,c}=-\half{}[[b^t,a],c]$ if $a,b,c\in\ppn$. 
\begin{example}\label{ex:surr**} 
We  will write out everything for the  case  $\gR=\fsu(r,r)$.
We began this example in \S\ref{ss2_surr} and continued
it in Example \ref{ex:surr*}. 

Now $\ppn$ identifies, in the obvious way, with the complex Jordan
algebra $M(r,\C)$ of
$r\times r$ matrices with Jordan product $A\circ B=\half(AB+BA)$
where $AB$ is the ordinary matrix product.  
The Jordan triple product is 
$\jtp{A,B,C}=\half(ABC+CBA)$.  We have
$\tr(A)=\Tr(A)$ and $F(A)=\Det(A)$ where $\Tr$ and $\Det$ are the
usual matrix trace and determinant. Also
$\pair{\left(\Atop{0}{0}\Atop{B}{0}\right),
\left(\Atop{0}{C}\Atop{0}{0}\right)}=\Tr(BC)$ and
$\left(\Atop{0}{0}\Atop{B}{0}\right)^t=
\left(\Atop{0}{B}\Atop{0}{0}\right)$.

To actually  to write  out our calculations
for this case with matrices, it is convenient to use the  basis 
$\{E_{i,j}\}$ of $M(r,\C)$ by elementary matrices (even though it is not
orthonormal). For instance, (\ref{eq:delvF}) becomes the familiar formula
$\pd{|z|}{z_{ij}}=|z|(z\inv)_{ji}$ where $|z|=\Det z$.
\end{example} 

\subsection{Computing  \boldmath $[\pila^y,[\pila^y,w]]$}  
\label{ss6_com} 
We now compute a bracket relation  in $\DtZ$.
\begin{lem}\label{lem:one} 
Let $y\in\pn$. Then  
\begin{equation}\label{eq:one} 
[\pila^y,w]=-w\del^y-2m(\la-\lao)[\del^yw]
\end{equation}
\end{lem}
\begin{proof}
The commutator  $[\pila^y,w]$ is a differential
operator on $\tZ$ of order at most $1$ and so we can        
write  it uniquely as the sum of  a vector field $\xi$ and a
function $g$. It is convenient to compute these parts 
individually. Using  (\ref{eq:pilaxxx*}) and (\ref{eq:delvw}) we find
\begin{equation}\nonumber
\begin{array}{lll}
\xi&=&-\textstyle\sum_{i,j}
\tr(\jtp{v_j,y,v_i}q)\,
\left([\del^{v_i}w]\del^{v_j}+[\del^{v_j}w]\del^{v_i}\right)\\[8pt]
&=&-w\textstyle\sum_{i,j}\tr(\jtp{v_j,y,v_i}q)
\tr(v_i\,q\inv)\del^{v_j}\\[8pt]
&=&-w\textstyle\sum_{j}\tr(\jtp{v_j,y,q\inv}q)\del^{v_j}\\[8pt]
&=&-w\textstyle\sum_{j}\tr(v_jy)\del^{v_j} 
=-w\del^y
\end{array}
\end{equation}
The fourth equality follows from Jordan identities. Indeed,
(i) the operator $\cP_{a,c}$ defined by $\cP_{a,c}(b)=\jtp{a,b,c}$  is
self-adjoint   and (ii) $\jtp{a,b,b\inv}=a$. Hence 
$\tr(\jtp{v_j,y,q\inv}q)=\tr(y\jtp{v_j,q,q\inv})=\tr(yv_j)$.

Next using (\ref{eq:pilaxxx*}), (\ref{eq:deldelw})
and self-adjointness of $\cP_{a,c}$ we find
\begin{eqnarray*}  
g&=&-\left(\tsum_{i,j}\tr(\jtp{v_j,y,v_i}q)
[\del^{v_i}\del^{v_j}w]\right)-2\nr\la[\del^yw]\\[6pt]
&=&-\tfour w\tr(\jtp{q\inv,y,q\inv}q)+\half w\tsum_i
\tr\left(\jtp{\jtp{q\inv,v_i,q\inv},y,v_i}q\right)
-m\la\,w\tr(yq\inv)\\[8pt]
&=&-(\tfour+m\la) w\tr(yq\inv)+\half w\tsum_i
\tr\left(yq\inv v_i^2\right)\\[8pt]
&=&-(\tfour+m\la-\half m)w\tr(yq\inv)
\end{eqnarray*}
For the third equality we used the identity
$\jtp{\jtp{q\inv,v,q\inv},q,v}=q\inv v^2$, and for the fourth we used
$\sum_iv_i^2=\nr\ove$  (see \cite[page 117]{F-K}).
\end{proof}

Notice that we can rewrite (\ref{eq:one})  as
\begin{equation}\label{eq:onep} 
[\pila^y,w]=-\del^yw-2m(\la-\laop)[\del^yw]
\end{equation}
 
\begin{lem}\label{lem:pi} 
Put $y=\ove_1$. Then  $[\pila^y,\del^y]= (\del^y)^2$.
\end{lem}
\begin{proof} 
Starting  from (\ref{eq:pilaxxx*})  we find
\begin{equation}\label{eq:[pi,de]} 
[\pila^y,\del^y]=\textstyle\sum_{i,j}
  \tr(\jtp{v_j,y,v_i}y)\del^{v_i}\del^{v_j}
 =\textstyle\sum_{i,j}
  \tr(yv_i)\tr(yv_j)\del^{v_j}\del^{v_i} = (\del^y)^2 
\end{equation}
We will explain the second equality.
We start from the   fact that $y=\ove_1$ is a primitive idempotent 
in the Jordan algebra $\pn$. Indeed,  $e=\sum_{i=1}^re_i$ is a
decomposition of $e$ into orthogonal primitive idempotents.

For any   primitive idempotent $y$, then  
the map $x\mapsto\jtp{y,x,y}$ is the orthogonal projection onto    
$\C{y}$ and so $\jtp{y,x,y}=y\tr(xy)$. Then
\begin{equation}\nonumber
y\tr(\jtp{v_j,y,v_i}y)=\jtp{y,\jtp{v_j,y,v_i},y} 
=\jtp{\jtp{y,v_j,y},v_i,y}
=y\tr(yv_j)\tr(yv_i) 
\end{equation}
because of the Jordan identity
$\jtp{a,\jtp{b,a,c},a}=\jtp{\jtp{a,b,a},c,a}$
 (see \cite[Ex. 8, page 40]{F-K}). This proves the second equality
in (\ref{eq:[pi,de]}).
\end{proof}
 
\begin{lem}\label{lem:two} 
Put    $y=\ove_1$. Then
\begin{equation}\label{eq:two} 
[\pila^y,[\pila^y,w]]=-m^2(\la-\lao)(\la-\laop)\,w\tr(yq\inv)^2
\end{equation}
Hence $[\pila^y,[\pila^y,w]]=0$ if and only if $\la$ equals $\lao$ or $\laop$.
\end{lem}
\begin{proof}
If $\la=\lao$ then  (\ref{eq:one}) gives
$[\pila^y,w]=-w\del^y$ and using Lemma \ref{lem:pi} we find
$[\pila^y,[\pila^y,w]]= w(\del^y)^2 -(w\del^y)\del^y=0$. 
If  $\la=\laop$ then (\ref{eq:onep})
gives $[\pila^y,w]=-\del^yw$ and we find
$[\pila^y,[\pila^y,w]]=-(\del^y)^2w+\del^y(\del^yw)=0$.
Thus (\ref{eq:two}) is true when $\la$ equals $\lao$ or $\laop$.  

Now we can prove (\ref{eq:two}) for all values of $\la$ without further
calculation by simply examining the form of $[\pila^y,[\pila^y,w]]$.
Let  $\la$ be arbitrary. Then   $\pila^y=\pilao^y-2\cla\del^y$ where
$\cla=m(\la-\lao)$.  Since $[\pilao^y,[\pilao^y,w]]=0$  we get
\begin{equation}\label{eq:four} 
[\pila^y,[\pila^y,w]]=-2\cla[\pi_{\lao}^y,[\del^yw]]
-2\cla[\del^y,[\pi_{\lao}^y,w]]+4\cla^2[\del^y[\del^yw]]
\end{equation}
The RHS is the sum of a vector field which in linear in $\la$ and a function
which is quadratic in $\la$. But the RHS vanishes for the two distinct values
$\lao$ and $\laop$. So the vector field   vanishes identically 
and the function is of the form
$(\la-\lao)(\la-\laop)g$  where  $g\in R(\tZ)$. Comparing  coefficients of
$\la^2$ we find $g=4m^2[\del^y,[\del^y,w]]$.
Using (\ref{eq:delvw}), (\ref{eq:delutr})   and the fact
that $y$ is a  primitive idempotent we find 
$4[\del^y,[\del^y,w]]=-w\tr(yq\inv)^2$.
\end{proof}

Lemmas \ref{lem:two} and  \ref{lem:wgfin} give (\ref{eq:wiff}).
This concludes  the proof of Theorem \ref{thm:la}.
\begin{remark}\label{rem:lowest}  
We can also write down   $Q^-$-semi-invariant (lowest weight) vectors 
$T\in\Wla[\lao]$ and $T'\in\Wla[\laop]$ using these methods.
We find that $T=w\delF$ and $T'=\delF w$.
This agrees with  \S\ref{ss4_sl2} since there $\lao=\four$
and $w\del_z=\half\del_w$.
\end{remark}

\subsection{Comparison of critical values}  
\label{ss6_comp}  
While it was easy to see that
$(\Elao,\tgalao,\pilao)$ and $(\E^{\laop},\tga_{\laop},\pi_{\laop})$
are anti-isomorphic (see Remark \ref{rem:thmla}), a more subtle fact is that
they are isomorphic.
We define an \emph{isomorphism of noncommutative models} from
$(\Ela,\tgala,\pila)$ to   $(\E^{\la'},\tga_{\la'},\pi_{\la'})$
to be  a filtered algebra  isomorphism $\sig:\Ela\to\E^{\la'}$
such that   $\sig(\pila^x)=\,\pi_{\la'}^x$ and  
$\gala=\ga_{\la'}(\gr\sig)$.
Then $\sig$ is $\g$-linear, and so
$G$-equivariant, with respect to the representations $\Pila$ and $\Pilap$.

\begin{prop}\label{prop:Innw} 
The inner automorphism $D\mapsto wDw\inv$ of $\DtZ$ maps $\Elaop$ onto
$\Elao$.  The induced map $\Innw:\E^{\laop}\to\Elao$
is an $\fS$-invariant   isomorphism  of noncommutative   models  from
$(\E^{\laop},\tga_{\laop},\pi_{\laop})$ to $(\Elao,\tgalao,\pilao)$.
\end{prop}    
\begin{proof} 
The result is clear once we prove  that 
$w\,\pi_{\laop}^xw\inv=\pi_{\lao}^x$ where $x\in\g$.
It suffices to check this for  $x\in\pp$ and $x\in\pn$.
Clearly $\Innw$ is the identity on $R(\tZ)$. So for $x\in\pp$ we find
$w\pi_{\laop}^xw\inv=w\pi^xw\inv=\pi^x=\pi_{\lao}^x$.
Next suppose $y\in\pn$. Then
\begin{equation}\nonumber 
w\,\pi^y_{\laop}w\inv=
\pi^y_{\laop}-[\pi^y_{\laop},w]w\inv 
=\pi^y_{\laop}+\del^y 
=\pi^y_{\laop-\frac{1}{2}{\nr}}=\pi^y_{\lao}  
\end{equation}  
The second equality follows by  (\ref{eq:one}) because 
$c_{\laop}=\half$  and so $[\pi^y_{\laop},w]=-\del^yw$. 
\end{proof}

We will find a sort of explanation  for   $\Innw$  later in 
\S\ref{ss7_comp}. Notice that $\Innw$ induces a filtered anti-isomorphism
$\Blaop\to\Blao$ which is then the restriction of the \emph{outer}
automorphism $D\mapsto wDw\inv$ of $\DZ$. It would be interesting
to give a direct geometric description of the 
corresponding anti-isomorphism $\Dlaop(X)\to\Dla(X)$ obtained by Fourier
transform.

\section{The noncommutative model    $\Elao$}  
\label{sec7}

\subsection{Algebraic structure of \boldmath $\Elao$}  
\label{ss7_alg}   
To begin with, we have
\begin{cor}\label{cor:Ela0=simple} 
$\Elao$ and $\Blao$ are   simple rings.
\end{cor}
\begin{proof} 
$\Dlao(X)$ is simple by  Proposition \ref{prop:Dlax=simple}   since   
$2\lao=1-\half[m]$   does  not lie in $\Z-\{1\}$. (Note $m\ge 1$ by Table
\ref{table:O}.) The result follows  by  Corollary \ref{cor:Ela=simple}.
\end{proof}

\begin{cor}\label{cor:Jlao} 
We have $\Jlao=\tau(\Jlao)=\Jlaop$. 
Moreover $\Jlao$ is a maximal  $2$-sided ideal in $\Ug$   
and its  infinitesimal character  is  given by the weights
$(-m\pm\half)\nu+\rho$.
\end{cor}
\begin{proof} 
$\Jlao=\Jlaop$ follows by Proposition \ref{prop:Innw}
while  $\Jlaop=\tau(\Jlao)$ follows by Corollary \ref{cor:tau_theta}.
$\Jlao$ is maximal since $\Blao\simeq\Ug/\Jlao$ is simple.
The infinitesimal character  follows by  Corollary \ref{cor:infch};
the two weights are then Weyl group conjugate.
\end{proof}
\begin{remark}\label{rem:McG}  
We checked that our infinitesimal character   is the
same as the one given by McGovern in \cite[Tables 5-10]{McG4}.
Moreover if the ``root multiplicity" $d$ given by 
$n=r+\binom{r}{2}d$ is equal to $2$, which happens exactly when
$\gR=\fsu(r,r)$, then our infinitesimal character coincides with $\half\rho$.
\end{remark}

\begin{cor}\label{cor:ann}  
Viewed as  $\Ug$-bimodules, both $\Blao$  and $\Elaood$  have left
annihilator and right annihilator  equal to $\Jlao$.
\end{cor}
\begin{proof} 
Both $\Blao$ and $\Elaood$ are faithful as right or left modules over
$\Blao$ since the algebra $\Elao$ has no zero-divisors. Since
$\Blao$ identifies with $\Ug/\Jlao$, the left and right annihilators in
$\Ug$ are $\Jlao$ and  $\tau(\Jlao)=\Jlao$.
\end{proof}

\subsection{A simple module for  \boldmath $\Elao$}  
\label{ss7_mod}   
Our construction of $\Elao$ gives us a natural module for it, namely  
$\RtZ$.  Our next result  produces a simple submodule.

Let $\cH$ be the $\Elao$-submodule of $\RtZ$ generated by the function
$1$.  

\begin{prop}\label{prop:cH=fs} 
$\cH$ is   a faithful simple $\Elao$-module. We have
\begin{equation}\label{eq:cH=} 
\cH=\Spp\oplus w\Spp
\end{equation}
The maximal Poisson abelian subalgebra $\cG$
\textup{(}from Corollary  \textup{\ref{cor:cG}}\textup{)} identifies with
$\cH$ under the   restriction homomorphism $\RM\to\RtZ$.
\end{prop} 

The proof requires two lemmas.

\begin{lem}\label{lem:cH=} 
$\cH$ is the subalgebra of $\RtZ$ generated by $\Spp$ and $w$.
The action of $\fS$ on $\RtZ$ induces an algebra  $\Z_2$-grading
$\cH=\cHev\oplus\cHod$   where $\cHev=\Spp$ and  $\cHod=w\Spp$
\textup{(}recall $w^2=F$\textup{)}.
\end{lem}
\begin{proof} 
Let $\A$ be the algebra generated by $\Spp$ and $w$.
Then $\A=\Spp\oplus w\Spp$ and this is the $\fS$-grading. 
Now $\A$, regarded as a space of
multiplication operators, lies in $\Elao_0$. Consequently $\A$, regarded    
as a space of functions, lies in $\cH$.
The problem then is to show that $\A$ is $\Elao$-stable.
We know (Corollary \ref{cor:Ela}) that
$\Elao$ is generated by $\Blao$ and $w$,
and $\Blao$ is generated by  $\pilao^x$, $x\in\g$.
The multiplication operator $w$ certainly preserves $\A$, and the
operators $\pilao^x$, $x\in\g$, preserve $\Spp$. So we need to check that
the operators  $\pilao^x$   preserve  $w\Spp$. For $x\in\q^+$, this is
clear. For $x\in\pn$,  this follows   using the bracket relations  
$[\pilao^x,w]=-w\del^x$ from  Lemma  \ref{lem:one}.
\end{proof}

\begin{remark}\label{rem:cH_stuff}  
Lemma \ref{lem:cH=} implies    that $\cH$ is the full subalgebra 
of all order zero differential operators (multiplication operators) in
$\Elao$, i.e., $\cH=\RtZ\cap\Elao$.
This suggests that  $\cH$  might be maximal abelian in $\Elao$;
we prove this in  Corollary \ref{cor:cG_cH}.
\end{remark}

Now  $\cHev$ and $\cHod$ are the $\Blao$-modules generated by $1$ and
$w$ respectively. Let $\tK$ be the double cover of $K$ which admits
$\sqrt{\chi}$ as a character.  
\begin{lem}\label{lem:cHevodd=low} 
$\cHev$ and $\cHod$ are
each lowest weight $\g$-representations and
generalized Verma modules for $\q^-$. 
The lowest weight vectors  
$1\in\cHev$ and $w\in\cHod$ are $\tK$-semi-invariant of weights
$\chi^{m-\half}$ and  $\chi^{m+\half}$ respectively.

$\cHev$ and $\cHod$ are  irreducible $(\g,\tK)$-modules with 
the same  annihilator $\Jlao$ in $\Ug$.  
\end{lem}
\begin{proof} 
By Corollary \ref{cor:verma} we know that $\cHev=\Spp$ is a 
Verma module for $\q^-$ with lowest weight vector $1$ of weight
$2m\lao\nu=(m-\half)\nu$.
Similarly  $\cHod=w\Spp$ is a Verma
module for $\q^-$ with lowest weight vector $w$ of 
weight $2m\laop\nu=(m+\half)\nu$. 
This follows because  for $x\in\kk$ we have 
$\pilao^x(w)=\nu(x)+2m\lao\nu(x)=(m+\half)\nu(x)$, and  
for $y\in\pn$ we have (by Lemma \ref{lem:one}) 
$\pilao^y(w)=([\pilao^y,w]+w\pilao)(1)=0$. 

The weights $(m\pm\half)\nu$ exponentiate to characters of $\tK$
and  then $\cHev$ and $\cHod$ are    $(\g,\tK)$-modules.
A theorem of Wallach (\cite{Wa}) says in particular that
$\Ug\otimes_{\U(\q^-)}\C_{s\nu}$ is irreducible as a $\g$-module 
if $s>m-1$.  Our values $s=m\pm\half$ satisfy this bound.
\end{proof}
 
\begin{remark}\label{rem:unit}  
The same theorem of Wallach says that $\cHev$ and $\cHod$ are
unitarizable as  $\gR$-representations.
We can  regard them as quantizations of the real nilpotent orbit $\OR$.
\end{remark}

\begin{proof}\textit{of Proposition \textup{\ref{prop:cH=fs}}.} 
$\cHev$ and $\cHod$ are  faithful simple  $\Blao$-modules by
Lemma \ref{lem:cHevodd=low}.
Since $\cHev$ and $\cHod$ carry different $\g$-representations, they 
are the only non-trivial $\Blao$-submodules of $\cH$. Neither           
$\cHev$ nor $\cHod$ is $\Elao$-stable, since multiplication by $w$ moves 
each into the other.  Thus $\cH$ is simple for $\Elao$. Faithfulness is
automatic as $\Elao$ is a simple ring.
Our descriptions of $\cG$ and $\cH$ in Corollary \ref{cor:cG} and Lemma
\ref{lem:cH=} imply that $\cG$ maps isomorphically onto $\cH$.
\end{proof}

\subsection{Using  \boldmath $\cH$ to realize  $\Elao$}  
\label{ss7_DD(cH)}   
Since $\Elao$ acts faithfully on $\cH$  we
have an algebra embedding
\begin{equation}\label{eq:Elao_into_End} 
\Elao\subset\End\sgfin(\cH)
\end{equation}
where the representation of $\g$ on $\End\cH$ is still given by
the operators $\Pilao^x$.
The algebra $\End\sgfin(\cH)$ is   much larger that $\Elao$;
in particular it contains $\End\sgfin(\cHev)\oplus\End\sgfin(\cHod)$.  
We next observe that  $\Elao$ is simply the subalgebra of
$\End\sgfin(\cH)$ consisting of differential operators.
\begin{cor}\label{cor:Elao==} 
$\Elao$ is the $\g$-finite part of $\DD(\cH)$ for the representation 
$\Pilao$.  
\end{cor}
\begin{proof} 
Clearly $\Elao$ lies in $\DD(\cH)\gfin$. The converse follows since 
$\RtZ$ is a localization of $\cH$ and so
any differential operator on   $\cH$ extends to one on $\RtZ$.
\end{proof}

We can also recover $\Elao$   as a vector space in the following  way.  
\begin{cor}\label{cor:Elao=vs} 
The natural map $\Elao\to\Hom\sgfin(\cHev,\cH)$
is a vector space isomorphism.
\end{cor}
\begin{proof} 
The map is injective because any differential operator on $\tZ$ is uniquely
determined by its  values on $\Spp$. 
This is true since any vector space basis of
$\pp$ is  a set of local \'etale coordinates on $\tZ$.
To prove surjectivity we need to show that 
if $L\in\Hom\sgfin(\cHev,\cH)$ then $L$ extends to  a differential
operator  $P$ on $\tZ$. 

We may write $L=L\ev+L\od$ where
$L\evod\in\Hom\sgfin(\cHev,\cH\evod)$.
(Read $\updownarrow$ like $\pm$.) Since  $L$ is $\g$-finite,
its components $L\ev$ and $L\od$ are each $\g$-finite,
and so in particular  are  $\pp$-finite. Now $x\in\pp$ acts by commutator 
with multiplication by $x$, i.e., $\Pilao^x(D)=[\pilao^x,D]=[x,D]$. 
It follows that  $L\ev$ lies in $\Dpn$. But also  
$\Pilao^x(w\inv L\od)=w\inv\Pilao^x(L\od)$ and so $w\inv L\od$ is
$\pp$-finite and thus lies in  $\Dpn$.  Let $P_1$ and $P_2$ be the 
differential operators on $Z$ defined by restriction of $L\ev$ and $w\inv
L\od$ respectively. Then $P=P_1+wP_2$ is the operator we wanted.
\end{proof} 

\subsection{Comparing  \boldmath $\Elao$ with $\Elaop$}  
\label{ss7_comp}   
Similarly, we can let   $\cHp\subset\RtZ$ be the $\E_{\laop}$-submodule  
generated by   $1$. Then $\cHp=\cHpev\oplus\cHpod$ where 
$\cHpev=\Spp$ and $\cHpod=w\inv\Spp $. 
Then  $1\in\cHpev$ and $w\inv\in\cHpod$ are
$\tK$-semi-invariant lowest weight vectors of weights
$\chi^{m+\half}$ and  $\chi^{m-\half}$ respectively.
(Notice that $\cHp$ is not  a subalgebra of $\RtZ$ and 
$\cHp$ does not lie in  $\E_{\laop}$.)
 
Now we get a nice way to derive the isomorphism $\Innw$ found in
Proposition \ref{prop:Innw}. Indeed we have an   isomorphism of
$(\g,\tK)$-modules
\begin{equation}\label{eq:} 
\cH'\to\cH,\quad f\mapsto wf
\end{equation}
which  carries $\cHev$ to $\cHpod$   and  $\cHod$ to $\cHpev$.  
Now $\Innw$ is simply the induced isomorphism 
$\End\sgfin(\cH')\to\End\sgfin(\cH)$ and this sends $\Elaop$ and $\Elao$.
 
\subsection{\boldmath The  algebra anti-automorphism $\be$} 
\label{ss7_be}

Let $\be$ be the composition  
\begin{equation}\label{eq:be} 
\Elao\,\xrightarrow{\delta}\,\E^{\laop}\,\xrightarrow{\Innw}\,\Elao
\end{equation}
Corollary \ref{cor:tdelta} and  Proposition \ref{prop:Innw} give
\begin{cor}\label{cor:be} 
$\be$ is an $\fS$-invariant anti-automorphism of   the
noncommutative model $(\Elao,\tgalao,\pilao)$. 
In particular $\be$ is $G$-invariant we have two commutative squares
\begin{equation}\label{eq:twosquares} 
\begin{array}{ccc}
\Ug&\mmaprightd{\pilao}&\Elao\\[8pt]
\mapdown{\tau}&&\mapdown{\beta}\\[8pt]
\Ug&\mmaprightd{\pilao}&\Elao 
\end{array}  \qquad \qquad
\begin{array}{ccc}
\gr\Elao&\mmaprightd{\gr\beta}&\gr\Elao\\[8pt]
\mapdown{\tgalao}&&\mapdown{\tgalao}\\[8pt]
\RM&\mmaprightd{\al}&\RM
\end{array}   
\end{equation}
$\be$  has order  $2$ or $4$; in fact $\be^2=1$ if $r$ is even while
$\be^2=\vsig$ if $r$ is odd where $\vsig$ is the
non-trivial element of $\fS$.
\end{cor}
\begin{proof} 
The first part is clear. Now $\be^2$ is a filtered $G$-invariant algebra
automorphism of $\Elao$ which is trivial on $\Blao$. 
By considering the induced action of $\be^2$ on
$\gr\Elao\simeq\RM$, we see easily that $\be^2$ lies in $\fS$. 
We have  $\be(w)=i^rw$ and so $\be^2=\vsig^r$.
\end{proof}
\begin{remark}\label{rem:be}  
In the sense of  (\ref{eq:twosquares}), (i)
$\be$ is the   unique  algebra anti-automorphism of $\Elao$ which
extends    $\tau$ , and (ii) 
$\beta$ is the unique $\g$-linear  endomorphism of $\Elao$ such
that $\gr\be$  corresponds to $\al$. 
\end{remark}
 
\begin{cor}\label{cor:be_descr} 
Here is a simple description of $\be$:
if   $\cL$ is  a subspace carrying an irreducible $G$-representation and
$\cL\subset\Elao_j$ with  $j$ as small as possible, then      
$\be(D)=i^{2j}D$  for all $D\in\cL$.
\end{cor}
This description is nice but does not reveal why  $\be$ is an 
anti-automorphism.
\subsection{\boldmath $\Elao$ is a superalgebra with supertrace}  
\label{ss7_T}   
The constant functions in $\Elao$ are the only $G$-invariants.
(This is true for a noncommutative model of any coadjoint orbit cover.)
So there is a unique $G$-linear map  
\begin{equation}\label{eq:7_T} 
T:\Elao\to\C  
\end{equation}
such that $T(1)=1$. The restriction of $T$ to $\Blao$ is a trace, i.e.,
$T(ab)=T(ba)$.  This is immediate  since $\Blao$ is a quotient of $\Ug$. 
Indeed $T(\pilao^xb-b\pilao^x)=0$ by  $\g$-invariance  and so 
$T(\pilao^{x_1}\cdots\pilao^{x_n}b)=T(b\pilao^{x_1}\cdots\pilao^{x_n})$.

We can ask now if $T$ is a trace on $\Elao$. The answer is no even for
$\g=\fsl(2,\C)$.   We will show this is remedied by introducing a
superalgebra structure on $\Elao$, which is \emph{filtered}  in the sense
that each filtration piece is the sum of its even and odd subspaces.
Corollary \ref{cor:be} gives
\begin{lem}\label{lem:super} 
Let $\E=\Elao$. Then $\E=\Eeven\oplus\Eodd$ where $\Eeven$ and
$\Eodd$ are the $\pm1$-eigenspaces of $\be^2$.   This makes $\E$ into a
$G$-equivariant filtered superalgebra.  

If  $r$ is odd then   $\Eeven=\E\ev$ and $\Eodd=\E\od$,
but if $r$ is even then $\Eeven=\E$ and $\Eodd=0$.
\end{lem}
An element $a\in\E$ is called
\emph{superhomogeneous} if $a$ lies in $\Eeven$ or $\Eodd$.
The \emph{parity} of $a$ is then $|a|=0$ or $|a|=1$ respectively.
\begin{prop}\label{prop:T} 
Our projection $T$ is a supertrace on $\Elao$, i.e.,
\begin{equation}\label{eq:Tab} 
T(ab)=(-1)^{|a||b|}T(ba)  
\end{equation}
when  $a$ and $b$ have same parity, while
$T(ab)=0$ when  $a$ and $b$ have different parity. 
\end{prop}
\begin{proof} 
Let $\E=\Elao$. 
Since $\E\simeq\RM$ is  multiplicity-free as a $G$-representation, 
there is a unique $G$-stable complement, call it $\E^j$,  
to   $\E_{j-\half}$ in  $\E_j$. Then $\E=\oplus_{j\in\halfN}\E^j$ is
a $\be$-stable vector space grading and   $\beta$
acts on $\E^j$ by multiplication by $i^{2j}$. So 
\begin{equation}\label{eq:Eevenodd} 
\Eeven=\oplus_{j\in\bbN}\,\E^j\AND
\Eodd=\oplus_{j\in\bbN+\half}\,\E^j
\end{equation}
 
The    pairing $\cP(a,b)=T(ab)$ is $G$-invariant. 
Suppose $j\neq k$. We know by Proposition \ref{prop:tO} that  
$\E^j$ and $\E^k$ contain no common  $G$-types 
and all $G$-types appearing are self-dual.  It follows by
Schur's Lemma that $\cP$ pairs  $\E^j$ and $\E^k$ trivially.
Now suppose $a,b\in\E^j$. Notice $T(\beta(c))=T(c)$ for any $c\in\E$. So 
\begin{equation}\label{eq:T(ab)=} 
T(ab)=T(\beta(ab))=T((\beta b)(\beta a))=(-1)^{2j}T(ba) 
\end{equation}
This proves $T$ is a supertrace.
\end{proof}

\begin{cor}\label{cor:cP} 
The  bilinear pairing $\cP(a,b)=T(ab)$ on $\Elao$
is $(\GfS)$-invariant, supersymmetric and non-degenerate.
\end{cor}
The pairing $\cP$ is \emph{supersymmetric} in the sense that
$\cP$ is symmetric on $\Deven$ and  is anti-symmetric on $\Dodd$,
while $\cP(a,b)=0$ if $a$ and $b$ have different parity.
Consequently, for any $\beta^2$-stable subspace $L$ in $\E$, 
the right and left $\cP$-orthogonal  subspaces coincide, thus
giving us a  notion of the  $\cP$-orthogonal subspace $L^\perp$. 
Now we say $\cP$ is \emph{non-degenerate} on   
$L$ if  $L^\perp\cap L=0$.
\begin{proof}
Everything is immediate except non-degeneracy.
Now $\E^\perp\cap\E=0$ because 
$\E$ is   simple   (Corollary \ref{cor:Ela0=simple})  and  
$\E^\perp\cap\E$ is a $2$-sided ideal in $\E$ 
which does not contain $1$.  
\end{proof}
\begin{remark}\label{rem:7_T}  
The multiplicity-free decomposition of $\E=\Elao$ into $G$-types is
orthogonal for $\cP$. Let $L$ be a $G$-type which lies in $\E^j$. If
$j\in\bbN$ then $\cP$ is symmetric non-degenerate on $L$ and so $L$ is an
orthogonal
$G$-representation; if If $j\in\bbN+\half$
then $\cP$ is anti-symmetric non-degenerate on $L$ and so $L$ is a
symplectic $G$-representation.  
\end{remark}

\subsection{\boldmath $\Ug$-bimodule structure of $\Elao$}  
\label{ss7_Ugbimod}    
As a first application of Proposition \ref{prop:T}, we get a quick proof of
another algebraic   fact about $\Elao$.

\begin{cor}\label{cor:bi_simple} 
$\Elaood$, like $\Elaoev=\Blao$,  is a simple  bimodule over $\Ug$. 
\end{cor}
\begin{proof}
We need to show that $\Elaood$ is a simple bimodule over 
$\Blao=\Ug/\Jlao$; we already know that $\Blao$ is a simple ring.
Suppose $L$ is a  $\Blao$-bisubmodule of $\Elaood$ and 
$a\in L$ with $a\neq 0$. Then $wa\in\Blao$. Now $\cP$ is non-degenerate
on $\Blao$ and so  there exists $b\in\Blao$ such that 
$\cP(b,wa)=T(bwa)=1$.
Then $\cP(ab,w)=T(abw)=-1$ since $T$ is a supertrace.
Thus $L$  is not  $\cP$-orthogonal to $\Wla$ (defined in \S\ref{ssE_gfin}).  
It follows,  since $\Wla$ is self-dual and appears only once in
$\Elaood$,  that $L$ contains $\Wla$. But $\Wla$ generates $\Elaod$ as a
bimodule over $\Blao$ by Corollary \ref{cor:Ela}(v). Thus  $L=\Elaood$.
\end{proof}

\section{Dixmier product on $\RM$}  
\label{sec8}
\subsection{Constructing the Dixmier  product}  
\label{ss8_circ}    
 
Suppose we have  a noncommutative model $(\E,\ga,\pi)$ of a graded
Poisson algebra $\RR$ with Hamiltonian symmetry $\g\to\RR[1]$,
$x\mapsto\phi^x$; see Definitions \ref{def:gPa} and \ref{def:ncm}. 
We say that $\bq:\RR\to\E$ is an associated 
\emph{quantization map} if (i) $\bq$ is $\g$-linear,   (ii) $\bq$ is filtered,
and (iii) the induced map $\gr\bq:\RR\to\gr\E$ is inverse to $\ga$. Here 
(i) means that $\bq(\{\phi^x,\psi\})=[\pi^x,\bq(\psi)]$ and  (ii) means that
$\bq(\RR[j])\subseteq\E_j$. So $\bq$ induces a vector space grading
$\E=\oplus_{j\in\halfN}\,\bq(\RR[j])$. In this way, we get a bijection 
between choices for $\bq$ and $\g$-linear gradings  
$\E=\oplus_{j\in\halfN}\E^j$ which satisfy
$\E_k=\oplus_{j\le k}\E^j$.

Our noncommutative model $(\Elao,\tgalao,\pilao)$ of $\RM$ admits a
unique   quantization map
\begin{equation}\label{eq:bq} 
\bq:\RM\to\Elao
\end{equation}  
because there is only one choice for the corresponding grading 
as $\Elao$ is multiplicity-free (cf. proof of Proposition \ref{prop:T}). 
We now get a new associative noncommutative product $\circ$  on $\RM$
defined by
\begin{equation}\label{eq:circ=} 
\phi\circ\psi=\bq\inv((\bq \phi)(\bq\psi))
\end{equation}
We say that $\circ$  is a \emph{Dixmier product} because it
makes $\RM$ into a Dixmier algebra for  $\M$.

\begin{example}\label{ex:sl2_bqcirc} 
In Example \ref{ss4_sl2}, $\bq$ is Weyl symmetrization map and  the circle
product   is   the Moyal star product specialized at $t=1$.
\end{example}
The  Euler grading defines a filtration of $\RM$  and also the projection 
\begin{equation}\label{eq:calT} 
\T:\RM\to\C
\end{equation}
Then $\T(\phi)$  is the   constant term of  $\phi$.  There is a 
supergrading   on $\RM$  given by
\begin{equation}\label{eq:RtOevenodd} 
\RMeven=\oplus_{j\in\bbN}\,\RM[j]\AND 
\RModd=\oplus_{j\in\bbN+\half}\,\RM[j]
\end{equation}
These are the  $\pm1$-eigenspaces of $\al^2$.

\subsection{Main theorem}  
\label{ss8_thm}    
We can now deduce
\begin{thm}\label{thm:circ} 
The Dixmier product $\circ$ is  $(\GcS)$-invariant and makes
$\RM$ into  a  filtered superalgebra 
where  \textup{(\ref{eq:RtOevenodd})} defines the supergrading.   With
respect to $\circ$,   $\al$ is an anti-automorphism   and   $\T$ is a
supertrace. The bilinear pairing
$\Q(\phi,\psi)=\T(\phi\circ\psi)$ on $\RM$ is
$(\GcS)$-invariant, supersymmetric, non-degenerate and  orthogonal for
the Euler grading.

Let $\RR[j]=\RM[j]$. Then, for all $j,k\in\halfN$,
\begin{equation}\label{eq:circ_jk} 
\RR[j]\circ\RR[k]\subseteq
\RR[j+k]\oplus\RR[j+k-1]\oplus\cdots\oplus\RR[|j-k|]
\end{equation}
Suppose $\phi\in\RR[j]$ and $\psi\in\RR[k]$ so that
$\phi\circ\psi=\sum_pC_p(\phi,\psi)$ where
$C_p(\phi,\psi)$ lies in $\RR[j+k-p]$. Then
\begin{eqnarray}
\phi\circ\psi&\equiv&\phi\psi+\half\{\phi,\psi\}\mod\RR[\le j+k-2]
\label{eq:circ_01}\\[3pt]
C_p(\phi,\psi)&=&(-1)^p\, C_p(\psi,\phi)
\label{eq:circ_par}
\end{eqnarray}
\end{thm}
\begin{proof} 
The map $\bq$ is $G$-invariant,   equivariant with
respect to  $\cS$ and $\fS$, and also $\bq$ intertwines $\al$ and $\be$.
So our results on $\Elao$ transfer over  to $\RM$ via $\bq$. This
proves the first paragraph. 

Since $\circ$ is a filtered superalgebra product, we have
$\RR[j]\circ\RR[k]\subseteq\bigoplus_{p\in\bbN}^{j+k}\RR[j+k-p]$. 
Now   proving  (\ref{eq:circ_jk}) reduces to showing that if
$\RR[j]\circ\RR[k]$ is not $\Q$-orthogonal to $\RR[s]$ then  $s\ge|j-k|$.
Showing this is easy since the hypothesis means that  
there exist $a\in\RR^j$, $b\in\RR^k$ and $c\in\RR^s$ such that
$T(abc)=1$. Then $bc$ has a component in $\RR^j$ and so
$k+s\ge j$. But also $T(bca)=\pm1$ (since $T$ is a supertrace) and so
similarly $s+j\ge k$. Hence $s\ge|j-k|$.

Since $\al$ is an anti-automorphism we find
\begin{equation}\label{eq:al_par} 
\al(\phi\circ\psi)=(\al\psi)\circ(\al\phi)=i^{2j+2k}\,\psi\circ\phi
\end{equation}
Then $i^{-2p}C_p(\phi,\psi)=C_p(\psi,\phi)$ and this proves 
(\ref{eq:circ_par}).  Next the relations  $C_0(\phi,\psi)=\phi\psi$ and  
$C_1(\phi,\psi)-C_1(\psi,\phi)=\{\phi,\psi\}$  follow
since $\gr\bq:\RR\to\gr\E$ is inverse to $\tgalao$.  
But $C_1(\phi,\psi)=-C_1(\psi,\phi)$ and so
$C_1(\phi,\psi)=\half\{\phi,\psi\}$. This proves (\ref{eq:circ_par}).
\end{proof}

Notice that  (\ref{eq:circ_jk})  implies that  if $\phi\in\RR[j]$ and
$\psi\in\RR[k]$ then
\begin{equation}\label{eq:T=C2j} 
\T(\phi\circ\psi)=\delta_{jk}\,C_{2j}(\phi,\psi)
\end{equation}

\begin{remark}\label{rem:bq}  
For   $\phi\in\RR[j]$ and  $\psi\in\RR[k]$ we have 
\begin{equation}\label{eq:bqPBm} 
\bq(\{\phi,\psi\})\equiv(\bq\phi)(\bq\psi)-(\bq\psi)(\bq\phi)
\mod\E_{j+k-2}
\end{equation}
This means that   $\bq$  approximately satisfies the Dirac  rule that 
quantization converts   Poisson brackets   into   commutators, thus $\bq$
deserves to be called a ``quantization map". 
\end{remark}

\begin{cor}\label{cor:exact} 
If $\phi\in\RM[1]$, for instance if $\phi=\phi^x$ where $x\in\g$, then     
for all $\psi\in\RM$  we have 
$\{\phi,\psi\}=\phi\circ\psi-\psi\circ\phi$. 
\end{cor}
\begin{proof} 

Let $\psi\in\RM[k]$. Then  we have 
$\phi\circ\psi=\phi\psi+\half\{\phi,\psi\}+C_2(\phi,\psi)$ by
(\ref{eq:circ_jk}). Now the result is immediate because of  
(\ref{eq:circ_par}). 
\end{proof}

\subsection{Underlying star product}  
\label{ss8_star}    
The  properties of the Dixmier product given in Theorem \ref{thm:circ}
reveal an underlying   star product. 
Here we mean \emph{star product} in the usual sense, 
except that we drop the requirement of locality and work with
 regular functions (rather than   all $\C$-valued smooth functions).
This  point of view is   known in star product theory. See 
e.g. \cite{CG}, \cite{ABC}, \cite{AB:starmin} for this and also for what it
means for a star product to be \emph{graded}  or \emph{strongly
invariant}.

\begin{cor}\label{cor:star} 
The Dixmier product $\circ$ is the  specialization at  $t=1$ of a
unique graded strongly $\g$-invariant star product $\star$ on $\RM$.
\end{cor}
\begin{proof} 
The graded star product is defined  by 
$\phi\star\psi=\sum_{p\in\bbN}C_p(\phi,\psi)t^p$ where
$\phi$ and $\psi$ are Euler homogeneous.
This is strongly $\g$-invariant by Corollary \ref{cor:exact}.
\end{proof}

\begin{remark}\label{rem:pos}
In \cite{me:2kpos}, we lift  the complex  conjugation
automorphism $\sig$ of $\OO$  (induced by a Cartan involution  of $\g$
which exchanges $\pp$ and $\pn$) to an antiholomorphic automorphism 
$\tsig$ of $\M$ (of order $2$ or $4$ according to whether $r$ is even or
odd) such that    (i) $\tsig$ induces a $\C$-antilinear $\circ$-algebra
automorphism of $\RM$ and (ii) the  pairing 
$(\phi|\psi)=\bbT(\phi\circ\psi^{\tsig})$ is  Hermitian positive-definite.
(In fact $(\cdot|\cdot)$ is Hermitian precisely 
because $\bbT$is a supertrace.)
Then $(\cdot|\cdot)$ is invariant under  the Lie algebra
$\{(x,x^\sig)\,|\,x\in\g\}$ and $\RM$ becomes a unitary representation    
of  $G$.
\end{remark}

\subsection{The operators \boldmath $\La^x$}   
\label{ss8_Lax} 
A natural first step in understanding the Dixmier product (or the
corresponding star product) is to compute the products $\phi^x\circ\psi$
where  $x\in\g$.  We get a neat form for the answer because of our
additional structure given by the supertrace etc.

Since our pairing $\Q$ on $\RM$ is supersymmetric and non-degenerate, it
makes sense to talk about the $\Q$-adjoint  of a linear endomorphism of
$\RM$. Theorem \ref{thm:circ} gives
\begin{cor}\label{cor:Lax} 
Let $\La^x$ be the $\Q$-adjoint of ordinary multiplication by $\phi^x$
where  $x\in\g$.  Then for every $\psi\in\RM$ we have
\begin{equation}\label{eq:phix*} 
\phi^x\circ\psi=\phi^x\psi+\half\{\phi^x,\psi\}+\La^x(\psi)
\end{equation}
The linear operators $\La^x$ satisfy:  
\begin{itemize}
\item[\rm(i)] $\La^x$ is graded of degree $-1$, i.e.,
$\La^x(\RM[j])\subseteq\RM[j-1]$.
\item[\rm(ii)] If $x\neq 0$ and $j$ is positive, then $\La^x$ is   non-zero
somewhere on $\RM[j]$.
\item[\rm(iii)] The operators $\La^x$ commute, i.e., 
$[\La^x,\La^y]=0$.
\item[\rm(iv)] The operators $\La^x$ transform in the
adjoint representation of $\g$, i.e., $[\Phi^x,\La^y]=\La^{[x,y]}$.
\item[\rm(v)] The operators $\La^x$ commute with the 
action of   $\cS$.
\end{itemize}
\end{cor} 
\begin{proof} 
Define $\La^x$ by $\La^x(\psi)=C_2(\phi^x,\psi)$.
Then  (\ref{eq:circ_jk}) and (\ref{eq:circ_01}) imply (\ref{eq:phix*}).
We claim $\Q(\phi^x\psi_1,\psi_2)=\Q(\psi_1,\La^x(\psi_2))$. 
We may assume $\psi_1\in\RR[j]$ and $\psi_2\in\RR[j+1]$. Then
(\ref{eq:phix*}) gives
$\Q(\phi^x\psi_1,\psi_2)=\Q(\phi^x\circ\psi_1,\psi_2)$  since the grading
of $\RR$ is $\Q$-orthogonal.   Now 
\[\Q(\phi^x\circ\psi_1,\psi_2)=T(\phi^x\circ\psi_1\circ\psi_2)
=T(\psi_1\circ\psi_2\circ\phi^x)=\Q(\psi_1,\psi_2\circ\phi^x)\]
since $T$ is a supertrace and $\phi^x$  and  $\psi_1\circ\psi_2$ are  even.
But $\Q(\psi_1,\psi_2\circ\phi^x)=\Q(\psi_1,C_2(\psi_2,\phi^x))$
and $C_2(\psi_2,\phi^x)=\La^x(\psi_2)$ by the parity relation
(\ref{eq:circ_par}).  This proves our claim.

Now the properties (i)-(v) of $\La^x$ follow  immediately from the
corresponding properties of their  $\Q$-adjoints, the  operators
$\psi\mapsto\phi^x\psi$. This works because of the properties of  $\Q$. 
\end{proof}
 
In \cite{me:2kpos} we give a formula for the operators $\La^x$. 

\begin{remark}\label{rem:Lax}  
We conjecture that there exist commuting homogeneous degree $-1$
algebraic differential operators $D^x$ on $\OO$ and a diagonalizable
$(\GcS)$-invariant algebraic differential operator $L$ on $\M$ with
positive real spectrum such that  
$\La^x=L\inv D^x$ as operators on $\RM$.

In the simplest case, where $\g=\fsl(2,\C)$, the operators $\La^x$ are 
differential. Indeed using the results in Examples \ref{ss4_sl2} and
\ref{ex:sl2_bqcirc} we find that 
$\Q(\xi^p,\ze^q)=\delta_{pq}\,2^{-p}p!$.
The operators $\La^x$ corresponding to the functions $\ze^2$, $\ze\xi$
and $\xi^2$ are 
$\four\pdb[2]{\xi^2}$, $-\four\pdb[2]{\xi\del\ze}$ and
$\four\pdb[2]{\ze^2}$.  
\end{remark}

If we identify $\RM$ with $\Elao$ via $\bq$, then the 
representation  (\ref{eq:Pi}) becomes 
\begin{equation}\label{eq:gog_RtO} 
\Pi:\gog\to\End_{\cS}\RM,\qquad 
\Pi^{(x,y)}(\psi)=\phi^x\circ\psi-\psi\circ\phi^y
\end{equation} 
Then $\Pi^{(x,x)}=\eta^x$ and Corollary \ref{cor:Lax} gives
\begin{cor}\label{cor:phi+La} 
For $x\in\g$ we have $\Pi^{(x,-x)}=2\phi^x+2\La^x$.
\end{cor}

\subsection{\boldmath  Dixmier product collapses on $\cG$}
\label{ss8_cG}   
Recall the maximal Poisson abelian
subalgebra $\cG$ of $\RM$ from Corollary \textup{\ref{cor:cG}}.
\begin{prop}\label{prop:*_cG} 
If $\psi$ and $\psi'$ lie in $\cG$  then
$\psi\circ\psi'=\psi\psi'$.
\end{prop}
\begin{proof} 
We can easily compute the restriction to $\cG$ of 
$\bq:\RM\to\Elao$.  We find
\begin{equation}\label{eq:bq_cG} 
\bq\,(\phi^P\ze^b)=Pw^b
\end{equation}
where $P\in\Spp$, $b\in\bbN$ and $\ze$ was defined in (\ref{eq:ze=}). 
Now   if  $\psi=\phi^P\xi^b$ and $\psi'=\phi^{P'}\xi^{b'}$ then
$\psi\circ\psi'=\ze\inv(Pw^bP'w^{b'})=\ze\inv(PP'w^{b+b'})
=\phi^{PP'}\xi^{b+b'}=\psi\psi'$.
\end{proof}

\begin{cor}\label{cor:cG_max_circ} 
$\cG$ is a maximal $\circ$-abelian subalgebra of $\RM$.
\end{cor}
\begin{proof} 
We have $\psi\circ\psi'=\psi\psi'$ and so $\cG$ is $\circ$-abelian.
Suppose $\phi\in\RM$ and $\phi\circ\psi=\psi\circ\phi$
for all $\psi\in\cG$. We can write $\phi=\sum_{j=0}^p\phi_j$
where $\phi_j\in\RM[j]$ and $\phi_p\neq 0$. 
Since $\cG$ is graded it follows
easily that $\{\phi_p,\psi\}=0$ for all $\psi\in\cG$.
But then $\phi_p\in\cG$ since $\cG$ is maximal Poisson abelian.
It follows by induction on $p$ that $\phi\in\cG$.
\end{proof}
\begin{cor}\label{cor:cG_cH} 
We have $\cH=\bq(\cG)$ and so $\cH$ is a maximal abelian subalgebra of 
$\Elao$.
\end{cor}

\bibliographystyle{plain}

\end{document}